\documentclass[10pt]{article}

\usepackage{amssymb}
\usepackage{amsmath}
\usepackage{enumerate}
\usepackage{latexsym}
\usepackage{mathrsfs}
\usepackage{amsthm}
\usepackage{verbatim}
\usepackage{graphicx}

\usepackage{float}

\usepackage[dvipsnames]{xcolor}
\usepackage{titlesec}
\usepackage{bm}
\usepackage[T1]{fontenc}

\usepackage{charter}

\usepackage{fancyhdr}
\usepackage{amscd}

\usepackage{subfig} 
\usepackage{caption}
\usepackage{algorithm} 
\usepackage{algpseudocode}%

\usepackage{tikz}
\usetikzlibrary{arrows.meta}
\usetikzlibrary{calc}

\usepackage[colorlinks=true,
linkcolor=blue,
citecolor=red,
urlcolor=blue]{hyperref}
\usepackage{cleveref}

\newtheorem{theorem}{Theorem}[section]

\newtheorem{lemma}[theorem]{Lemma}
\newtheorem{proposition}[theorem]{Proposition}

\newtheorem{remark}[theorem]{Remark}
\newtheorem{definition}[theorem]{Definition}

\numberwithin{equation}{section}
\renewcommand\theequation{\oldstylenums{\thesection}%
	.\oldstylenums{\arabic{equation}}}

\newcommand{\nullspace}{\text{null}}

\usepackage{pgfplots}
\pgfplotsset{compat=1.16}

\newcommand{\R}{\mathbb R}
\newcommand{\Pk}{P_{k}}
\newcommand{\Pkp}{\Pk^+}
\newcommand{\Pkm}{\Pk^-}
\newcommand{\Pkpm}{\Pk^\pm}

\newcommand{\edit}[1]{{\color{black} #1}}

\usepackage{geometry}
\newcommand{\email}[1]{\href{mailto:#1}{#1}}
\crefname{remark}{remark}{remarks}

\title{Constrained optimization and higher-order\\rigidity of bar frameworks}
\author{Xuenan Li\thanks{Department of Mathematics, Purdue University (\email{xuenan@purdue.edu}).}
\qquad Christian D. Santangelo\thanks{Department of Physics, Syracuse University (\email{cdsantan@syr.edu}).}
\qquad Miranda Holmes-Cerfon\thanks{Department of Mathematics, University of British Columbia (\email{holmescerfon@math.ubc.ca}).}}
\hypersetup{pdftitle={Constrained optimization and higher-order rigidity of bar frameworks},pdfauthor={Xuenan Li, Christian D. Santangelo, Miranda Holmes-Cerfon}}
\date{}

\begin{document}
\maketitle
\begin{abstract}
We establish a connection between constrained optimization problems, and conditions for higher-order rigidity of bar frameworks. 
We then use this connection to design an optimization algorithm to construct \edit{isostatic and under-constrained} bar frameworks that are rigid but not first-order rigid. Specifically,  prestress stable 
frameworks arise as the solution to a simple optimization problem, which asks to maximize or minimize the length of one edge while keep the other edge lengths fixed. 
By starting with a random first-order rigid framework, we can thus design a wide variety of prestress stable frameworks. 
We then present a bifurcation method to design frameworks that are third-order rigid. 
Our results highlight connections between concepts in rigidity theory and constrained optimization, offering new insights into the construction and analysis of bar frameworks with higher-order rigidity.

\end{abstract}
\begin{quote}\small

\end{quote}

\section{Introduction}\label{sec:intro}

Bar frameworks are graphs consisting of rigid bars (edges) connected by flexible hinges (vertices). These models serve as powerful abstractions for analyzing the structural and dynamic behavior of interconnected systems across various scales. Their applications range from understanding the micro and macro structures of mechanical metamaterials \cite{hutchinson2006structural, li2023some,li2025effective,li2025nonlinear} to designing metamaterials with unique properties \cite{bertoldi2017flexible, borcea2017new, paulose2015selective, rayneau2018analytic} to exploring the mathematical principles of origami folding \cite{filipov2017bar, demaine2007geometric} to investigating jammed particle arrangements \cite{connelly2019rigidity, henkes2016rigid, liu2010jamming}, to understanding allosteric mechanisms in biology \cite{rocks2017designing, yan2017architecture}, to studying molecular properties \cite{holmes2016enumerating, sartbaeva2006flexibility, whiteley2005counting}, and to analyzing protein structures \cite{jacobs2001protein}.

A central question across these areas is whether a given framework is locally rigid, meaning that any continuous motion of the vertices that preserves edge lengths is a rigid-body motion. 
The mathematical theory of rigidity \cite{connelly1980rigidity,connelly1996second} offers several sufficient conditions for local rigidity. The strongest condition is first-order rigidity \cite{asimow1978rigidity}. Mechanically, a framework which is first-order rigid resists forces in a manner that is ``typical'' -- a small deformation leads to changes at first order in the lengths of the bars. A weaker condition is prestress stability. Mechanically, this corresponds to the framework exhibiting first-order behaviour when a particular set of tensions is put on its edges. 
There are even weaker notions of rigidity, but they all share the property that small deformations can cause the bar lengths to change at higher order in the size of the deformation. Mechanically, such frameworks feel shaky. 
We will call any framework which is rigid but not first-order rigid \emph{higher-order rigid}.

Higher-order rigid frameworks are of interest, partly as illustrations of mathematical principles, but also because of their practical or scientific applications. A prominent example is provided by tensegrity frameworks -- rigid assemblies of bars and cables stabilized by self-stress -- which are closely tied to second-order rigidity and enable lightweight, adaptive designs such as NASA's Superball Bot \cite{tibert2011review,connelly1996second,shah2022tensegrity,garanger2021soft,vespignani2018design}. They also arise in observations of clusters of small particles, where the additional entropy associated with the higher-order motions can stabilize such clusters thermodynamically \cite{meng2010free,kallus2017free}. 
Another example arises in the so-called ``vertex model'', which considers two-dimensional graphs with additional constraints on the areas of the cells spanned by vertices. Second-order rigid configurations in the vertex model  have been used to understand the transition between fluid-like and solid-like behaviour in living tissues \cite{aspinwall2025rigidity,bi2015density,merkel2019minimal}.

\begin{figure}[!htb]
	\centering
	
	\subfloat[]{
		\begin{tikzpicture}[scale=0.6]
			\coordinate (A2) at (0,0);
			\coordinate (B2) at (-0.5, 0.87);
			\coordinate (C2) at (0, 1.75);
			\coordinate (D2) at (1, 1.75);
			\coordinate (E2) at (1.5, 0.87);
			\coordinate (F2) at (1, 0);
			
			\draw[ultra thick] (A2) -- (B2); 
			\draw[ultra thick] (B2) -- (C2); 
			\draw[ultra thick] (C2) -- (D2); 
			\draw[ultra thick] (D2) -- (E2); 
			\draw[ultra thick] (E2) -- (F2); 
			\draw[ultra thick] (A2) -- (F2); 
			\draw[ultra thick] (B2) -- (E2); 
			\draw[ultra thick] (C2) -- (F2); 
			\draw[ultra thick] (A2) -- (D2); 
			
			\fill[black] (A2) circle (3.5pt);
			\fill[black] (F2) circle (3.5pt);
			\fill[black] (B2) circle (3.5pt);
			\fill[black] (C2) circle (3.5pt);
			\fill[black] (D2) circle (3.5pt);
			\fill[black] (E2) circle (3.5pt);
		\end{tikzpicture}
	}
	\hfil
	\subfloat[]{
		\begin{tikzpicture}[scale=0.6]
			\coordinate (A2) at (0, 0);
			\coordinate (B2) at (0, 1.8);
			\coordinate (C2) at (0.8, 7.8484e-01);
			
			\coordinate (D2) at (2.4, 0);
			\coordinate (E2) at (2.4, 1.8);
			\coordinate (F2) at (1.6, 7.8484e-01);
			
			\draw[ultra thick] (A2) -- (B2); 
			\draw[ultra thick] (B2) -- (C2); 
			\draw[ultra thick] (A2) -- (C2); 
			
			\draw[ultra thick] (D2) -- (E2); 
			\draw[ultra thick] (E2) -- (F2); 
			\draw[ultra thick] (D2) -- (F2); 
			
			\draw[ultra thick] (B2) -- (E2); 
			\draw[ultra thick] (A2) -- (D2); 
			\draw[ultra thick] (C2) -- (F2); 
			
			\fill[black] (A2) circle (3.5pt);
			\fill[black] (D2) circle (3.5pt);
			
			\fill[black] (B2) circle (3.5pt);
			\fill[black] (C2) circle (3.5pt);
			\fill[black] (E2) circle (3.5pt);
			\fill[black] (F2) circle (3.5pt);
		\end{tikzpicture}
	}\hfil
	\subfloat[]{
		\begin{tikzpicture}[scale=0.6]
			\coordinate (A1) at (0,0);
			\coordinate (A2) at (1,1);
			\coordinate (A3) at (0,2);
			
			\coordinate (B1) at (1.5,0);
			\coordinate (B2) at (0.5,1);
			\coordinate (B3) at (1.5,2);
			
			\draw[ultra thick] (A1) -- (A2); 
			\draw[ultra thick] (A2) -- (A3); 
			\draw[ultra thick] (A3) -- (A1); 
			
			\draw[ultra thick] (B1) -- (B2); 
			\draw[ultra thick] (B2) -- (B3);
			\draw[ultra thick] (B3) -- (B1);
			
			\draw[ultra thick] (A2) -- (B2); 
			\draw[ultra thick] (A1) -- (B1);
			\draw[ultra thick] (A3) -- (B3);
			
			\fill[black] (A1) circle (3.5pt);
			\fill[black] (A2) circle (3.5pt);
			\fill[black] (A3) circle (3.5pt);
			
			\fill[black] (B1) circle (3.5pt);
			\fill[black] (B2) circle (3.5pt);
			\fill[black] (B3) circle (3.5pt);
		\end{tikzpicture}
	}\hfil
	\subfloat[]{
		\begin{tikzpicture}[scale=.35]
			\coordinate (A1) at (0, 0);
			\coordinate (A2) at (4,0);
			\coordinate (A3) at (4,4);
			\coordinate (A4) at (0,4);
			
			\coordinate (L1) at (2,0);
			\coordinate (L2) at (4,2);
			\coordinate (L3) at (2,4);
			\coordinate (L4) at (0,2);
			
			\coordinate (B1) at (2, 1);
			\coordinate (B2) at (3, 2);
			\coordinate (B3) at (2, 3);
			\coordinate (B4) at (1, 2);
			
			\draw[ultra thick,black] (A1) -- (A2); 
			\draw[ultra thick,black] (A2) -- (A3); 
			\draw[ultra thick,black] (A3) -- (A4); 
			\draw[ultra thick,black] (A4) -- (A1); 
			
			\draw[ultra thick,black] (B1) -- (B2); 
			\draw[ultra thick,black] (B2) -- (B3);
			\draw[ultra thick,black] (B3) -- (B4);
			\draw[ultra thick,black] (B4) -- (B1); 
			
			\draw[ultra thick] (L1) -- (B1); 
			\draw[ultra thick,black] (L2) -- (B2);
			\draw[ultra thick,black] (L3) -- (B3);
			\draw[ultra thick,black] (L4) -- (B4); 
			
			\fill[black] (A1) circle (6pt);
			\fill[black] (A2) circle (6pt);
			
			\fill[black] (A3) circle (6pt);
			\fill[black] (A4) circle (6pt);
			
			\fill[black] (L1) circle (6pt);
			\fill[black] (L2) circle (6pt);
			\fill[black] (L3) circle (6pt);
			\fill[black] (L4) circle (6pt);
			
			\fill[black] (B1) circle (6pt);
			\fill[black] (B2) circle (6pt);
			\fill[black] (B3) circle (6pt);
			\fill[black] (B4) circle (6pt);
	\end{tikzpicture}}
	\hfil
	\subfloat[]{
		\includegraphics[width=0.12\linewidth]{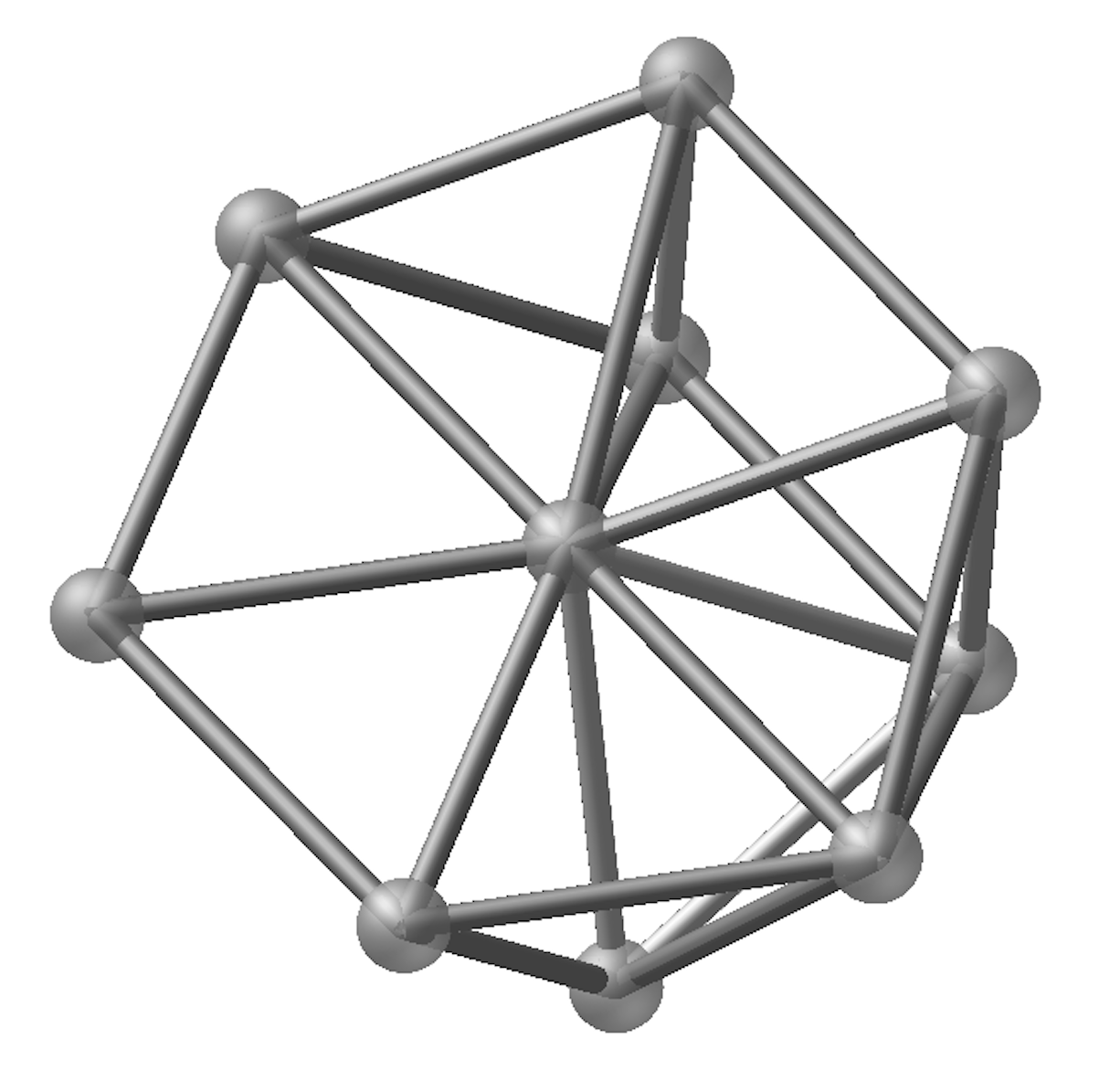}
	}
	\hfil
	\subfloat[]{
		\includegraphics[width=0.12\linewidth]{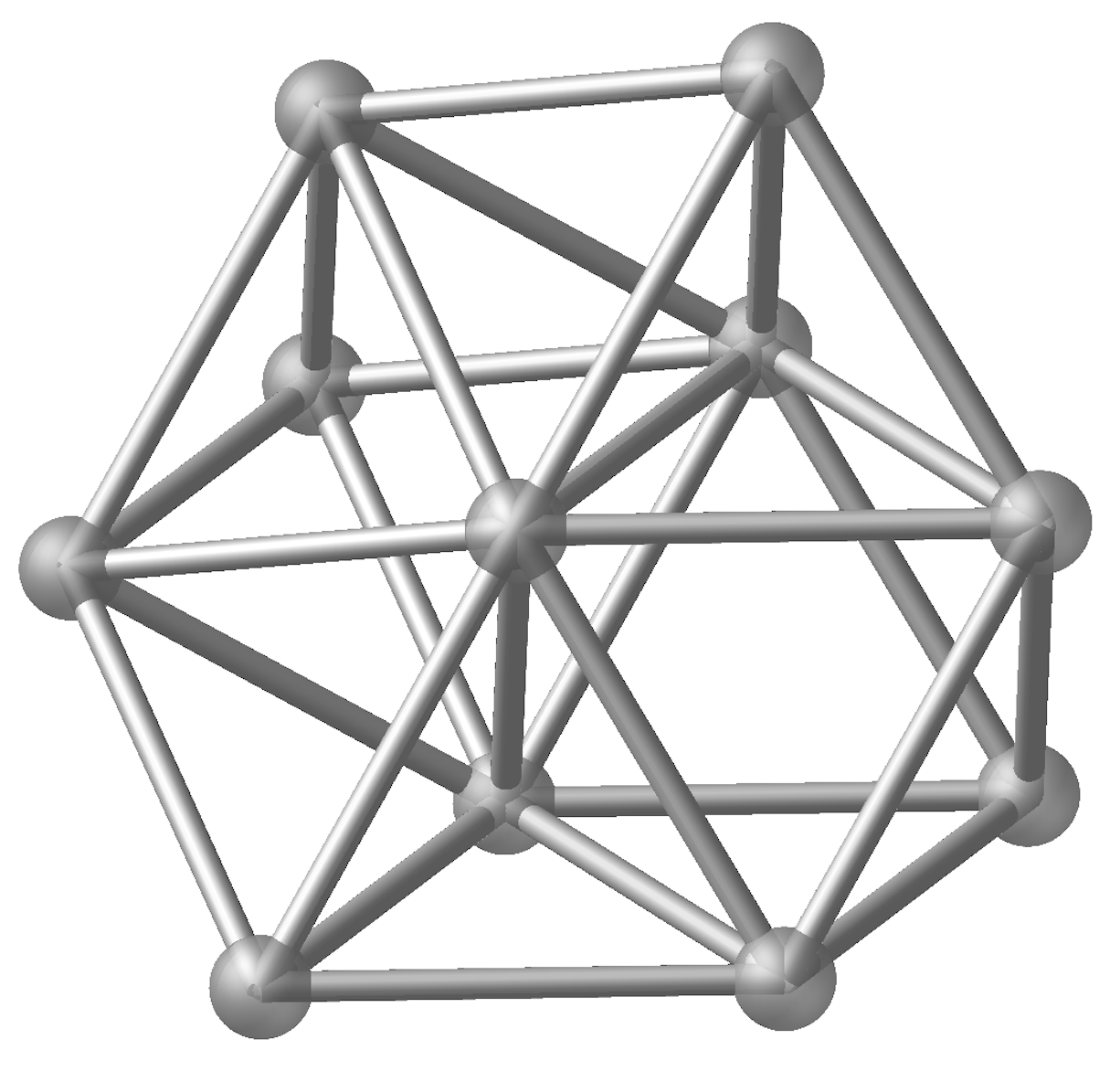}
	}\\
	\subfloat[]{
		\begin{tikzpicture}[scale=0.6]
			\coordinate (A2) at (0,0);
			\coordinate (B2) at (0.75*-6.1425e-01, 7.8911e-01);
			\coordinate (C2) at (0.75*4.8962e-01, 1.4758e+00);
			\coordinate (D2) at (0.75*2.4634e+00, 1.7987e+00);
			\coordinate (E2) at (0.75*3.3594e+00, 1.0005e+00);
			\coordinate (F2) at (0.75*1, 0);
			
			\draw[ultra thick] (A2) -- (B2); 
			\draw[ultra thick] (B2) -- (C2); 
			\draw[ultra thick] (C2) -- (D2); 
			\draw[ultra thick] (D2) -- (E2); 
			\draw[ultra thick] (E2) -- (F2); 
			\draw[ultra thick] (A2) -- (F2); 
			\draw[ultra thick] (B2) -- (E2); 
			\draw[ultra thick] (C2) -- (F2); 
			\draw[ultra thick] (A2) -- (D2); 
			
			\fill[black] (A2) circle (3.5pt);
			\fill[black] (F2) circle (3.5pt);
			\fill[black] (B2) circle (3.5pt);
			\fill[black] (C2) circle (3.5pt);
			\fill[black] (D2) circle (3.5pt);
			\fill[black] (E2) circle (3.5pt);
		\end{tikzpicture}
	}
	\hfil
	\subfloat[]{
		\begin{tikzpicture}[scale=0.6]
			\coordinate (A2) at (0, 0);
			\coordinate (B2) at (0.75*8.5812e-01, 9.7654e-01);
			\coordinate (C2) at (0.75*1.8396e+00, 7.8484e-01);
			
			\coordinate (D2) at (0.75*3.2000e+00, 0);
			\coordinate (E2) at (0.75*3.6146e+00, 2.1606e+00);
			\coordinate (F2) at (0.75*2.4417e+00, 9.3002e-01);
			
			\draw[ultra thick] (A2) -- (B2); 
			\draw[ultra thick] (B2) -- (C2); 
			\draw[ultra thick] (A2) -- (C2); 
			
			\draw[ultra thick] (D2) -- (E2); 
			\draw[ultra thick] (E2) -- (F2); 
			\draw[ultra thick] (D2) -- (F2); 
			
			\draw[ultra thick] (B2) -- (E2); 
			\draw[ultra thick] (A2) -- (D2); 
			\draw[ultra thick] (C2) -- (F2); 
			
			\fill[black] (A2) circle (3.5pt);
			\fill[black] (D2) circle (3.5pt);
			
			\fill[black] (B2) circle (3.5pt);
			\fill[black] (C2) circle (3.5pt);
			\fill[black] (E2) circle (3.5pt);
			\fill[black] (F2) circle (3.5pt);
		\end{tikzpicture}
	}\hfil
	\subfloat[]{
		\begin{tikzpicture}[scale=0.6]
			\coordinate (A1) at (0,0);
			\coordinate (A2) at (2.2087e-01,0.75*1.6856e+00);
			\coordinate (A3) at (-1.1883e+00, 0.75*2.1995e+00);
			
			\coordinate (B1) at (4.9815e-01, 0);
			\coordinate (B2) at (-4.2870e-01, 0.75*1.5430e+00);
			\coordinate (B3) at (9.8204e-01, 0.75*2.9607e+00);
			
			\draw[ultra thick] (A1) -- (A2); 
			\draw[ultra thick] (A2) -- (A3); 
			\draw[ultra thick] (A3) -- (A1); 
			
			\draw[ultra thick] (B1) -- (B2); 
			\draw[ultra thick] (B2) -- (B3);
			\draw[ultra thick] (B3) -- (B1);
			
			\draw[ultra thick] (A2) -- (B2); 
			\draw[ultra thick] (A1) -- (B1);
			\draw[ultra thick] (A3) -- (B3);
			
			\fill[black] (A1) circle (3.5pt);
			\fill[black] (A2) circle (3.5pt);
			\fill[black] (A3) circle (3.5pt);
			
			\fill[black] (B1) circle (3.5pt);
			\fill[black] (B2) circle (3.5pt);
			\fill[black] (B3) circle (3.5pt);
		\end{tikzpicture}
	}\hfil
	\subfloat[]{
		\begin{tikzpicture}[scale=.35]
			\coordinate (A1) at (0, 0);
			\coordinate (A2) at (4,0);
			\coordinate (A3) at (2.6004e+00,3.7471e+00);
			\coordinate (A4) at (-1.3996e+00,3.7471e+00);
			
			\coordinate (L1) at (2,0);
			\coordinate (L2) at (3.3002e+00,1.8736e+00);
			\coordinate (L3) at (6.0036e-01,3.7471e+00);
			\coordinate (L4) at (-6.9982e-01,1.8736e+00);
			
			\coordinate (B1) at (1.7710e+00, 8.7585e-01);
			\coordinate (B2) at (1.9975e+00, 1.7698e+00);
			\coordinate (B3) at (8.5574e-01, 3.2366e+00);
			\coordinate (B4) at (1.8475e-01, 1.7795e+00);
			
			\draw[ultra thick,black] (A1) -- (A2); 
			\draw[ultra thick,black] (A2) -- (A3); 
			\draw[ultra thick,black] (A3) -- (A4); 
			\draw[ultra thick,black] (A4) -- (A1); 
			
			\draw[ultra thick,black] (B1) -- (B2); 
			\draw[ultra thick,black] (B2) -- (B3);
			\draw[ultra thick,black] (B3) -- (B4);
			\draw[ultra thick,black] (B4) -- (B1); 
			
			\draw[ultra thick] (L1) -- (B1); 
			\draw[ultra thick,black] (L2) -- (B2);
			\draw[ultra thick,black] (L3) -- (B3);
			\draw[ultra thick,black] (L4) -- (B4); 
			
			\fill[black] (A1) circle (6pt);
			\fill[black] (A2) circle (6pt);
			
			\fill[black] (A3) circle (6pt);
			\fill[black] (A4) circle (6pt);
			
			\fill[black] (L1) circle (6pt);
			\fill[black] (L2) circle (6pt);
			\fill[black] (L3) circle (6pt);
			\fill[black] (L4) circle (6pt);
			
			\fill[black] (B1) circle (6pt);
			\fill[black] (B2) circle (6pt);
			\fill[black] (B3) circle (6pt);
			\fill[black] (B4) circle (6pt);
	\end{tikzpicture}}
	\hfil
	\subfloat[]{
		\includegraphics[width=0.12\linewidth]{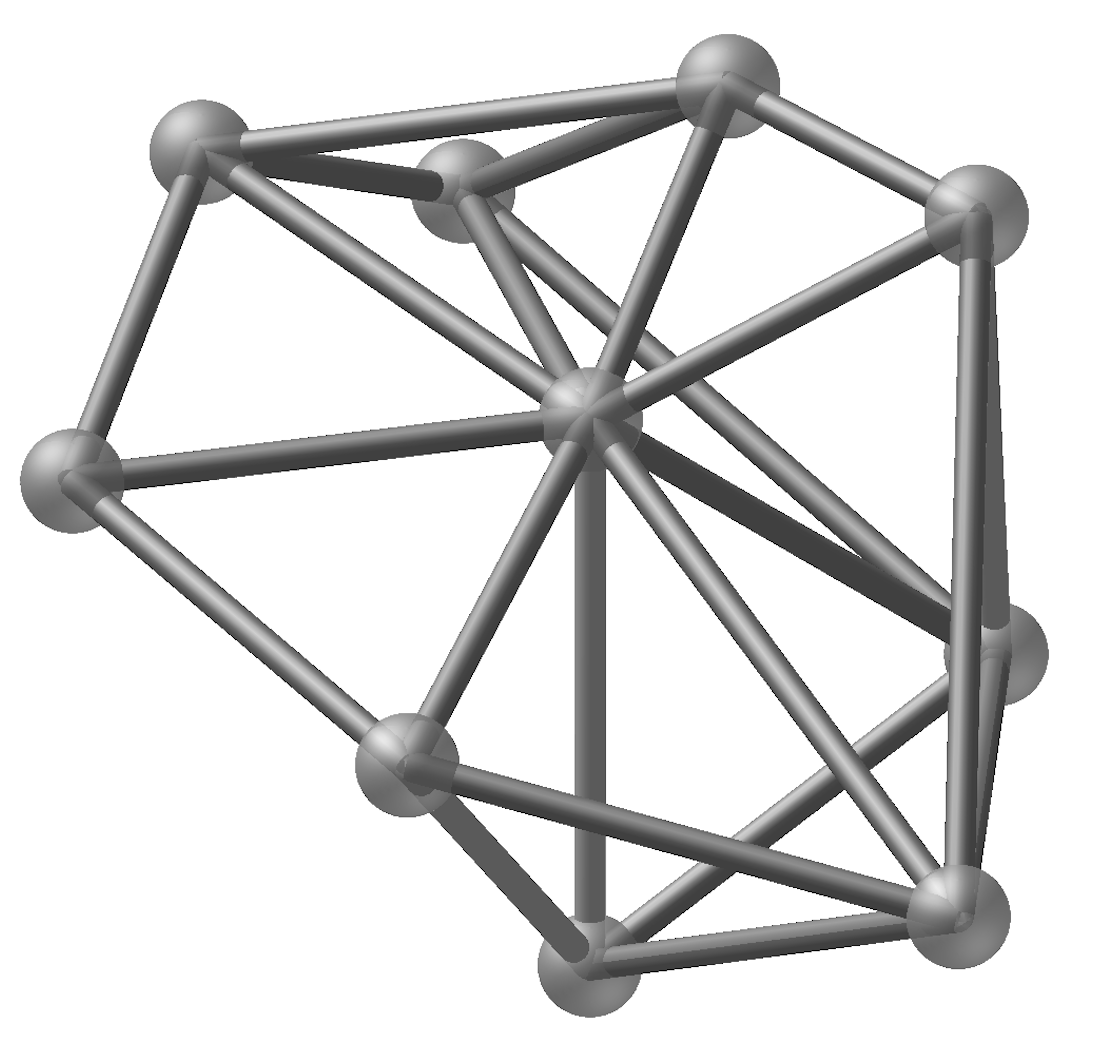}
	}
	\hfil
	\subfloat[]{
		\includegraphics[width=0.12\linewidth]{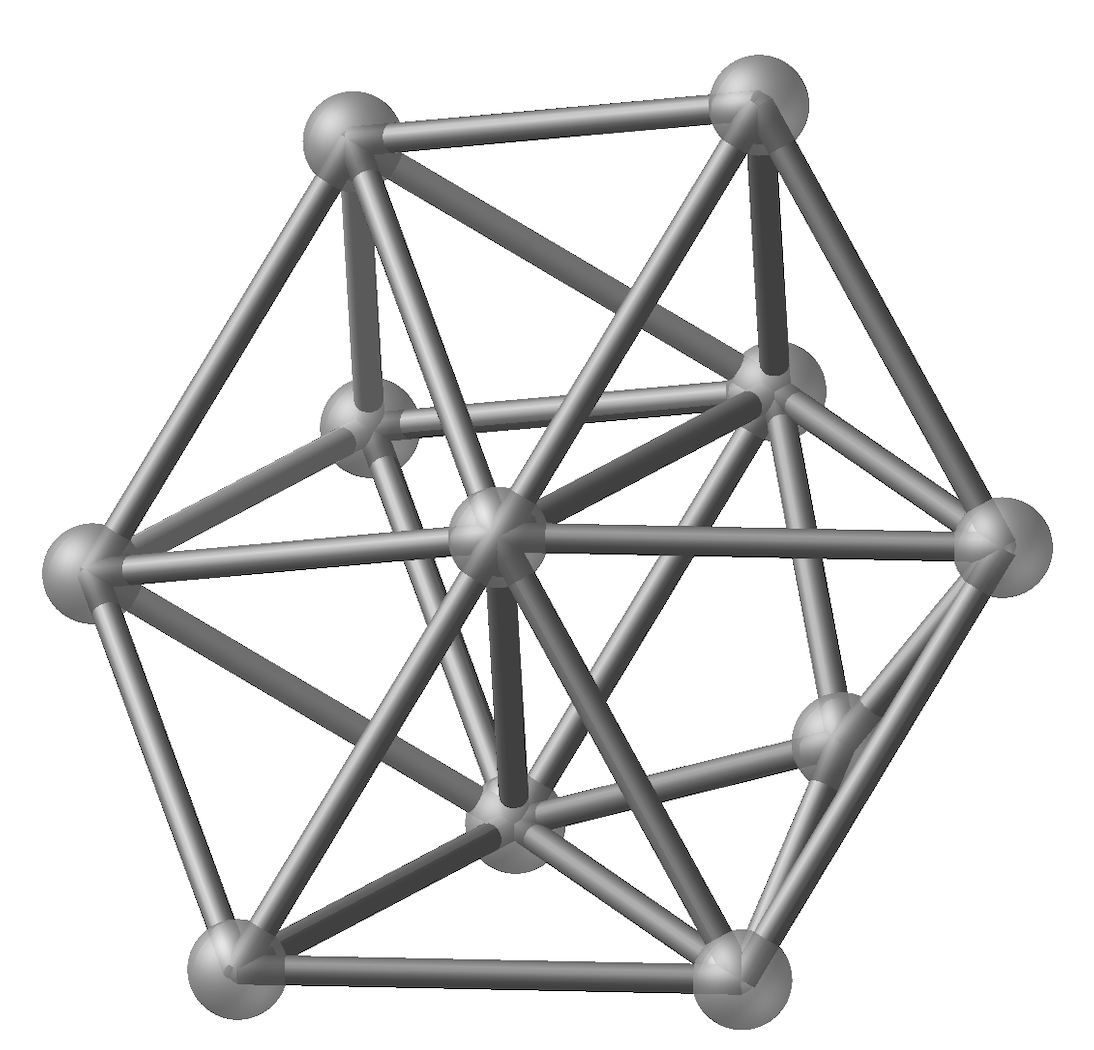}
	}
	\caption{Various prestress-stable bar frameworks that are not first-order rigid: the first row depicts well-known examples -- (a) from \cite{connelly1996second}, (b)-(c) from \cite{grasegger2019graphs}, (d) from \cite{roback2025tuning}, and (e)-(f) from the sphere-packing data in \cite{holmes2016enumerating}; the second row shows prestress-stable frameworks obtained by perturbing these well-known structures and then applying our constrained optimization approach.}
	\label{fig:zoo-rigid-examples}
\end{figure}
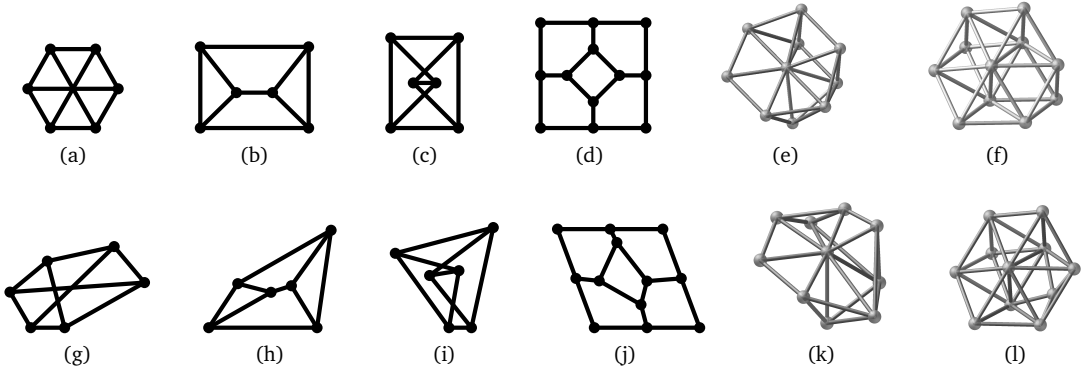

Yet, examples of higher-order rigid frameworks are not common. 
%
For a fixed graph, if edge lengths are chosen randomly, the framework is generally either first-order rigid or flexible, but not higher-order rigid. There are examples of higher-order rigid  frameworks known in the literature, however, these are often constructed by hand, and most rely on special symmetries. For example,  they might have parallel bars, as in the examples in  \cref{fig:zoo-rigid-examples}(a)-(d), or edges with the same lengths or coplanar vertices, as in \cref{fig:zoo-rigid-examples}(e) and (f). 

Our goal is to develop a systematic method to construct such frameworks, while also gaining insight into how ``hard'' it is to produce them.
We have found that it is in fact not hard to produce prestress stable frameworks -- that they arise robustly as the solution to a simple optimization problem. 
Let us explain this problem. Consider a $d$-dimensional  framework with $n$ vertices and $m$ edges. Let $p$ denote the collection of vertex positions, let $p_{i,1},p_{i,2}\in \R^d$ be the endpoints of the $i$th edge, and let $l_i$ be the length of the $i$th edge. Suppose we try to minimize the length of the first edge, while maintaining the lengths of the remaining edges:
\begin{equation}\label{eq:opt1}
	\min_{p \in \mathbb{R}^{nd}} \: f_1(p) = |p_{1,1} - p_{1,2}|^2, \quad
	\text{s.t } \quad f_i(p) = |p_{i,1} - p_{i,2}|^2 = l_i^2, \quad i = 2,\dots,m.
\end{equation}
Under mild conditions, a local solution $p^*$ corresponds to a framework that is prestress stable. \edit{For isostatic and under-constrained frameworks, such a solution also admits a nontrivial infinitesimal flex (Theorem~\ref{thm:prestress-stability}).} If instead we try to maximize the length of the first edge, we obtain another prestress stable framework. This leads to an easy method to create examples of prestress stable frameworks with (nearly) arbitrary edge lengths: start with a framework that is first-order rigid, remove an edge, and then use any constrained optimization procedure to find local maxima and local minima of the length of the removed edge. 

Some examples of prestress stable frameworks found by this method are shown in \cref{fig:zoo-rigid-examples}(g)-(l). Notably, these examples have edge lengths that are all distinct, and they lack any obvious symmetries.

\edit{Though constrained optimization has been used to study bar frameworks with interesting rigidity properties~\cite{gorkavyi2018birosettes,gorkavyy2016jessen,nawratil2022snappability}, our contribution is to establish a rigorous connection between rigidity and optimality for the constrained optimization problem~\eqref{eq:opt1}. In particular, we relate first-order optimality condition to self-stress and second-order optimality conditions to prestress stability. Specifically,} we show that the Karush-Kuhn-Tucker (KKT) conditions in constrained optimization, which provide a Lagrange multiplier $\lambda^*$ for the optimization problem \eqref{eq:opt1} and are closely linked to the concept of a critical point for the problem, are equivalent to the condition for a self-stress $\omega$ in rigidity theory (Lemma \ref{lem:KKT}). Therefore, under mild conditions we may take $\omega=\lambda^*$, so that the existence of a self-stress in rigidity theory, is equivalent to the KKT conditions for a constrained optimization problem. Another link arises at second-order: we show that the second-order sufficient conditions for the constrained optimization problem \eqref{eq:opt1}, are equivalent to the condition for prestress stability of the corresponding framework (Theorem \ref{thm:prestress-stability}). 
The latter is what allows us to prove prestress stability of the solution to \eqref{eq:opt1}.

Furthermore, we show that a converse of the optimality result also holds, and this leads to an interesting characterization of prestress stable frameworks: consider a prestress stable framework with self-stress $\omega$ that certifies prestress stability. Then, every edge $E_k$ for which $\omega_k>0$, is a local minimum for a problem of the form \eqref{eq:opt1} but with $|p_{k,1}-p_{k,2}|^2$ as the objective. Similarly, every edge $E_k$ for which $\omega_k<0$, is at a local maximum! Therefore, returning to the optimization problem \eqref{eq:opt1}, once we have optimized the length of one edge, all other edges bearing a self-stress are individually at a local optimum  (Theorem \ref{thm:2nd-sufficient-cond}). 

We extend these results to consider a broader class of higher-order rigid frameworks in several ways. One, we generalize the optimization problem \eqref{eq:opt1}, to allow prescribing components of the self-stress on certain edges. We show that one can either choose the length of an edge, or the self-stress on that edge. This generalizes the force-density method which is widely used in the engineering literature to build structures where all the components of a self-stress are prescribed (e.g. \cite{hain2025optimizing,schek1974force}). 
Second, we show how to design frameworks that are so-called ``third-order rigid'' \cite{gortler2025higher} by adjusting another edge length, to create a critical point of the length function $|p_{1,1}-p_{1,2}|^2$ which behaves locally as a cubic (Theorem \ref{thm:cubic-growth}). For this, we rely on bifurcation theory, which says that when a local maximum and a local minimum of a function merge, the function acquires cubic behaviour.  
Third, we explore a method to design frameworks containing multiple self-stresses, by finding the global optimum of a particular objective function created from the sub-determinants of the rigidity matrix. However, this method comes with fewer guarantees than our previous methods.

Our results are related to others that have appeared in the literature. A similar constrained optimization approach to constructing tensegrity structures was suggested in \cite{tibert2011review}, but without a rigorous guarantee of the prestress stability. 
A similar approach to ours was also followed in \cite{li2025constrained}, which showed that if instead of local optima one considers saddle points of $f_1(p)$ subject to the constraints, the resulting  framework is flexible with configuration spaces consisting of intersecting branches. 
Our idea was inspired by studies of the vertex model, used as a model for tissue dynamics, which stretch out the model until it reaches a critical point where the system becomes prestress stable \cite{bi2015density,merkel2019minimal}.
The idea that bifurcation points, or catastrophes, could be useful for studying frameworks or related systems has been suggested earlier; for example \cite{heaton2022catastrophe} developed a method based on numerical algebraic geometry to compute the set of catastrophe points, though the rigidity properties of these points have not been studied.

Our paper is organized as follows. In Section \ref{sec:prelim}, we review basic concepts in rigidity theory and in constrained optimization. In Section \ref{sec:our-theorems}, we introduce our constrained optimization approach within the context of bar frameworks, state and prove the main theorems. In Section \ref{sec:examples}, we present examples of prestress stable but not first-order bar frameworks obtained using our method. In Section \ref{sec:higher-rigidity}, we describe an approach based on bifurcation theory to design a third-order rigid framework.

\section{Background information on rigidity and optimization}\label{sec:prelim}

In this section, we review the key ideas from rigidity theory in \cref{subsec:prelim-rigidity} and constrained optimization in \cref{subsec:constrained-opt} that we will later connect. The rigidity results have been extensively detailed in \cite{asimow1978rigidity,connelly1980rigidity,connelly1996second} and our summary loosely follows that of \cite{holmes2021almost}. 
Our review of constrained optimization is based on \cite{nocedal1999numerical}.

\subsection{A brief review of  rigidity theory for bar frameworks}\label{subsec:prelim-rigidity}

A bar framework $G = (p,E)$ in $\mathbb{R}^d$ is a list of vertices $p$ embedded in $\R^d$, and a set of edges $E$. We write the vertices as a vector $p = (p_1, p_2, \dots, p_n) \in \mathbb{R}^{nd}$, where $p_i\in \R^d$ is the coordinates of vertex $i$, and $n$ is the total number of vertices. The edge set is $E = \{E_1,E_2,\dots,E_m\}$ where $m$ is the number of edges. 
We sometimes index the vertices based on their relationship to the edges: for each edge $E_i$ with $i=1,2,\dots,m$, we write the vertices it connects as $p_{i,1},p_{i,2}$.  

Each edge $E_i$ has a length $l_i=|p_{i,1}-p_{i,2}|$, so the vertices satisfy the system of algebraic equations 
\begin{equation}\label{eqn:bar-constraints}
	f_i(p)=l_i^2 \quad \text{with}\quad f_i(p):=\left|p_{i,1}-p_{i,2}\right|^{2}, \qquad i=1,2,\dots,m.
\end{equation}
For a given $p$, the framework can be translated,  rotated, or reflected without changing the lengths of the edges, so there is at minimum a $\frac{d(d+1)}{2}$-dimensional subspace of solutions to \eqref{eqn:bar-constraints}, corresponding to applying an isometry of $\R^d$ to the vertices. We say a $q$ is \emph{congruent to $p$} if it is obtained by applying such an isometry to $p$. Handling these isometries will present some technical though not fundamental challenges in our theory, though see \cite{gortler2025higher} for an in-depth study of how to handle them within the context of rigidity theory. 

If there are no other solutions to \eqref{eqn:bar-constraints} with the same edge lengths, then the framework is rigid. Specifically:

\begin{definition}
	A framework $G=(p,E)$ is \emph{(locally) rigid} if there is a neighborhood $U$ of $p$ such that any $q\in U$ with the same edge lengths as $p$, is congruent to $p$.
\end{definition}

Typically, one would expect a framework to be rigid when it has $m\geq nd-\frac{d(d+1)}{2}$ edges. Such an argument, based on equating the number of variables and the number of constraints, goes back at least to Maxwell \cite{maxwell1864calculation}, and has been used extensively in the study of disordered physical systems. Although this argument isn't true in general, it still gives a useful way of categorizing frameworks based on the number of edges, which we will frequently refer to.

\begin{definition}
	We say a framework with $n$ vertices and $m$ edges $E$ is \emph{isostatic} when $m + \frac{d(d+1)}{2} = nd$;  \emph{under-constrained} when $m + \frac{d(d+1)}{2} < nd$, and  \emph{over-constrained} when $m + \frac{d(d+1)}{2} > nd$. 
\end{definition}

We will mostly be interested in isostatic and under-constrained frameworks in this paper. 

Testing for rigidity is co-NP hard \cite{Abbot-Hard} and there is no known efficient algorithm for a general framework. Therefore, stronger flavors of rigidity have been developed. Of relevance for this paper will be first-order rigidity, the strongest notion, and prestress stability, a weaker notion that is still possible to test efficiently. 

To introduce these notions, consider the following system of equations, to be solved for a vector $v \in \mathbb{R}^{nd}$, indexed below using edges:
\begin{equation}\label{eqn:first-order-test}
	(p_{i,1}- p_{i,2}) \cdot (v_{i,1}- v_{i,2}) = 0, \qquad i = 1,2,\dots,m.
\end{equation}
These come from assuming an analytic deformation of the framework $p(t)$ with $p(0) = p$,  taking $\frac{d}{dt}$ of \eqref{eqn:bar-constraints} and evaluating at $t=0$. The vector $v = p'(0)$ has the interpretation of the velocity of the vertices required to deform the framework initially. The system \eqref{eqn:first-order-test} can be written in matrix-vector form as 
\begin{equation}\label{eqn:first-order-test-matrix}
	R(p) v = 0,
\end{equation}
where the matrix $R(p)\in \R^{m\times nd}$, called the \emph{rigidity matrix}, is defined so its $i$th row  is  $\nabla f_i(p)$, i.e.
\begin{equation}\label{R}
	R(p)^T = \begin{pmatrix}
		\nabla f_1(p)^T & \cdots & \nabla f_m(p)^T
	\end{pmatrix}.
\end{equation}
An element $v\in \nullspace R(p)$ is called an \emph{infinitesimal flex}.

It is worth noting that \eqref{eqn:first-order-test-matrix} admits a $\tfrac{d(d+1)}{2}$-dimensional space of solutions corresponding to infinitesimal isometries. We call such solutions the \emph{trivial flexes}, and denote the space of trivial flexes at $p$ as $\mathcal{T}(p)$ (sometimes abbreviated as $\mathcal T$). 
Any other solution $v$ which does not lie in $\mathcal T(p)$ is  a \emph{nontrivial infinitesimal flex}. 
We will sometimes need a subspace of such nontrivial infinitesimal flexes, so we will consider $v\in \nullspace R(p) \cap \mathcal T(p)^\perp$. (Note that in our theory that follows one can replace $\mathcal T^\perp$ with any space complementary to $\mathcal T$.) For brevity, we will call such a $v\in \nullspace R(p) \cap \mathcal T(p)^\perp$ a \emph{$\mathcal T^\perp$-flex}.

\begin{definition}
	A framework $G=(p,E)$ is \emph{first-order rigid} if it \edit{has no} nontrivial infinitesimal flex.
\end{definition}

There are many proofs that a first-order rigid framework is rigid (e.g. \cite{asimow1978rigidity}). First-order rigidity is the strongest notion of rigidity, and it is straightforward to test using linear algebra, however it is too strong for many examples of interest. 

A weaker notion of rigidity comes from looking at second-order deformations. The second-order equations, motivated by taking $\frac{d^2}{dt^2}\Big|_{t=0}$ of \eqref{eqn:bar-constraints}, are 
\begin{align}
	(p_{i,1}- p_{i,2}) \cdot (a_{i,1}- a_{i,2}) + (v_{i,1}- v_{i,2}) \cdot (v_{i,1}- v_{i,2}) = 0, \qquad i = 1,2,\dots,m, \label{eqn:second-order-test}
\end{align}
where $a \in \mathbb{R}^{m \times nd} = p''(0)$ can be interpreted as the initial acceleration of the vertices. The matrix-vector form of \eqref{eqn:second-order-test} is
\begin{align}
	R(p) a = -R(v) v. \label{eqn:second-order-matrix}
\end{align}
A framework is \emph{second-order rigid} if there is no solution $(v,a)$ to the system \eqref{eqn:first-order-test-matrix} and \eqref{eqn:second-order-matrix} with a nontrivial infinitesimal flex $v$. A second-order rigid  framework is rigid, though the proof is nontrivial \cite{connelly1980rigidity} (see also \cite{gortler2025higher}). 

In between first- and second-order rigidity is prestress stability. 

\begin{definition}
	A framework $G=(p,E)$ is \emph{prestress stable} if there exists a $w \in \nullspace \: R^T(p)$ such that 
	\begin{align}
		\omega^{\top} R(v) v>0 \qquad \text{for all } v \in \nullspace \: R(p) \cap \mathcal{T}(p)^\perp.\label{eqn:prestress-stable}
	\end{align}
\end{definition}
Briefly, to see where this definition comes from, consider the Fredholm alternative for \eqref{eqn:second-order-matrix}, which says there is no solution $a$ provided that, for each $v \in \nullspace \: R(p) \cap \mathcal{T}(p)^\perp$, there is a $\omega_v \in \nullspace \: R^T(p)$ such that $\omega_v^{\top} R(v) v>0$. Prestress stability asks for a single $\omega \in \nullspace \: R^T(p)$ that works for all $v$.

A framework that is prestress stable is second-order rigid, and hence rigid \cite{connelly1996second}. Prestress stability can be tested efficiently using semidefinite programming \cite{holmes2021almost}.

Vectors in the left null space of $R(p)$ play an important role, theoretically and physically. 

\begin{definition}
	A vector $\omega \in \nullspace \: R^T(p)$ is called a \emph{self-stress}. 
\end{definition}

A general vector $\omega \in \mathbb{R}^{m}$ has the interpretation of a \emph{stress}, which is an assignment of tensions to the edges. Component $w_i$ represents the tension on edge $E_i$ with the sign indicating whether the edge is under tension or compression. The net force on vertices is $R^T(p) w$. 
When $\omega$ is a self-stress, then there is no net force on the vertices, so the framework achieves local force equilibrium. 

\begin{remark}\label{rmk:fundamental-lin-alg}
	The number of independent self-stresses and of $\mathcal{T}^{\perp}$-flexes are related, by the fundamental theorem of linear algebra:
	\begin{equation}\label{eqn:fundamental-lin-alg}
		\mbox{dim}(\text{self-stresses}) - \mbox{dim}(\text{$\mathcal{T}^{\perp}$-flexes}) = nd-m-\frac{d(d+1)}{2}.
	\end{equation}
	Hence, for an isostatic framework, the number of self-stresses and the number of $\mathcal{T}^{\perp}$-flexes are equal. Typically, a ``random'' isostatic framework will have no $\mathcal{T}^{\perp}$-flexes and hence no self-stresses. A framework which is under-constrained will have $\mathcal{T}^{\perp}$-flexes, hence it is not first-order rigid, but only in special cases will it also have self-stresses. Hence, for isostatic or under-constrained frameworks, the presence of a self-stress is non-generic, and indicates that it could be rigid at higher order than first-order. 
\end{remark}

\subsection{A brief review of equality-constrained optimization}\label{subsec:constrained-opt}
For a set of edges $E$ and a corresponding list of lengths $\{l_i\}_{i\in E}$, consider the following constrained optimization problem: 
\begin{equation}\label{eqn:opt}
	(\Pkpm) \qquad  \min_{p \in \mathbb{R}^{nd}}  \pm f_k(p) \quad
	\text{s.t } \quad  f_i(p)=l_i^2, \quad i = 1,\dots,m, \quad i\neq k.
\end{equation}
The above notation refers concisely to two separate problems: problem $\Pkp$ has objective function  $f_k(p)$, and problem $\Pkm$ has objective function $-f_k(p)$. 
Recall $f_i(p) = |p_{i,1} - p_{i,2}|^2$ is the length of the $i$th edge.
This problem asks to minimize ($\Pkp$) or maximize ($\Pkm$) the length of the $k$th edge, subject to the constraints that the lengths of the remaining edges are fixed. A point $p$ which satisfies the constraints is called a \emph{feasible} point.

Recall that a solution to problem $\Pkpm$ is defined as follows. 

\begin{definition}\label{def:localmin}
	We say $p^*$ is a \emph{local solution  of $\Pkp$ (or $\Pkm$)} if (i) $p^*$ is feasible, and (ii) there exists a neighborhood $U$ of $p^*$ such that any point $q\in U$ satisfying the constraints has $f_k(q)\geq f_k(p^*)$ (or $f_k(q)\leq f_k(p^*)$).  
\end{definition}

Notice that because the rigid-body motions do not change the objective nor the constraints, a local solution will never be strict -- it will never be possible to replace the condition $f_k(q)\geq f_k(p^*)$ with $f_k(q)> f_k(p^*)$. Therefore, it is convenient to consider the solution modulo rigid-body motions, which we define as follows.

\begin{definition}\label{def:strictlocalmin}
	We say $p^*$ is a \emph{strict local solution  of $\Pkp$ ($\Pkm$) up to rigid body motions} if (i) $p^*$ is feasible, and (ii) there exists a neighborhood $U$ of $p^*$ such that any point $q\in U$ satisfying the constraints, and such that $q$ is not congruent to $p^*$, has $f_k(q)> f_k(p^*)$ ($f_k(q)< f_k(p^*)$). 
\end{definition}

The theory of constrained optimization gives various necessary or sufficient conditions for a point $p^*$ to be a local solution, strict or not, of $\Pkpm$. 
We survey the relevant results here, presented for the specific form of our equality constraints in $\Pkpm$ and modified slightly to account for the rigid-body motions. 


The key condition that most solutions must satisfy, and that will link optimality conditions to ideas from rigidity,   is the following. 

\begin{definition}\label{def:KKT}
	We say $p \in \mathbb{R}^{nd}$ satisfies the \emph{Karush-Kuhn-Tucker conditions} for $\Pkpm$, or the \emph{KKT conditions},  if (i) $p$ is feasible for $\Pkpm$, and (ii) there exists a vector 
	$\lambda  \in \mathbb{R}^{m}$, called the \emph{Lagrange multiplier}, with $\lambda_k=1$, and such that
	\begin{equation}\label{eqn:KKT}
		\sum_{i=1}^m \lambda_i \nabla f_i(p) = 0. 
	\end{equation}
\end{definition}

\begin{remark}
	Our definition is slightly different from how the KKT conditions are  usually presented, because we wish to be retain flexibility regarding which edge is being optimized without introducing cumbersome notation. The usual definition of the Lagrange multiplier, say for $k=1$, would have $\lambda=(\lambda_2,\ldots,\lambda_m)\in \R^{m-1}$ and would impose condition $\nabla f_1(p) + \sum_{i=2}^m\lambda_i \nabla f_i(p)=0.$
\end{remark}

The function 
\begin{equation}\label{L}
	L(p,\lambda) = \sum_{i=1}^m \lambda_i  f_i(p)
\end{equation}
is the \emph{Lagrangian} for the constrained optimization problem \eqref{eqn:opt}. The KKT conditions ask that $\nabla_p L= 0$, where $\nabla_p$ represents the gradient with respect to the $p$-variables. Hence, a point satisfying the KKT conditions can be thought of as a critical point for the function $f_k(p)$ subject to the constraints in \eqref{eqn:opt}.

When the constraints in \eqref{eqn:opt} meet the following regularity assumption, the KKT condition must hold at a local solution $p^*$ of \eqref{eqn:opt}.
\begin{definition} \label{def:LICQ}
	Given a point $p \in \mathbb{R}^{nd}$ that is feasible for $\Pkpm$, we say the \emph{linear independence constraint qualification (LICQ)} holds at $p$ if  the gradients of the constraints, $\{\nabla f_i(p)\}_{i\neq k}$, are linearly independent.
\end{definition}

\begin{theorem}[First-order necessary conditions; \cite{nocedal1999numerical}, Theorem 12.1]\label{thm:firstordernecessary}
	Suppose $p^*$ is a local solution of $\Pkpm$, and the LICQ holds. Then the KKT conditions hold for some Lagrange multiplier $\lambda^*$. 
\end{theorem}

\begin{remark}[Unique Lagrange multiplier]\label{rmk:unique-lagrange-multiplier}
	When the LICQ and KKT conditions hold at $p$, the Lagrange multiplier must be unique. To see why, suppose there are two nonequal Lagrange multipliers $\lambda,\gamma$. Then, subtracting conditions \eqref{eqn:KKT} gives $\sum_{i=1}^m (\lambda_i - \gamma_i) \nabla f_i(p) = 0$. Since $\lambda \neq \gamma$ but $\lambda_k=\gamma_k=1$, the vectors $\{\nabla f_i(p)\}_{i\neq k}$ are not linearly independent, contradicting LICQ.
\end{remark}

The KKT condition in Theorem \ref{thm:firstordernecessary} (which doesn't need to be modified to account for rigid-body motions) gives a necessary condition that a local solution must satisfy, but it is not sufficient. Indeed, it does not distinguish a local minimum or a local maximum, and it also holds for a saddle point for the constrained problem. To guarantee we have a local solution will require so-called second-order conditions. The conditions for a strict local solution must be modified from the standard conditions to account for rigid-body motions. 

\begin{definition}
	Given a point $p^*$ satisfying the KKT conditions \eqref{eqn:KKT} for $\Pkpm$, the \emph{critical cone up to rigid body motions} is
	\begin{equation}\label{eqn:cone-rigid}
		C_k(p^*) = \big\{v \in \mathcal T(p^*)^\perp\: |\: \nabla f_i(p^*) v = 0, \quad i\neq k\big\}.
	\end{equation}
\end{definition} 

\begin{theorem}[Second-order sufficient conditions]\label{thm:secondordersufficient}
	Consider problems $\Pkpm$ with feasible point $p^*$, and suppose there is a Lagrange multiplier $\lambda^*$ such that the KKT conditions \eqref{eqn:KKT} are satisfied. If
	\begin{equation}\label{eqn:secondsufficient}
		v^T \Big(\sum_{i=1}^m \lambda_i^* \nabla^2 f_i(p^*)\Big) v > 0 \qquad \forall v \in {C_k(p^*)}, \quad v\neq 0,
	\end{equation}
	then $p^*$ is a strict local solution up to rigid body motions for $\Pkp$.
	If instead the inequality is reversed,  then $p^*$ is a strict local solution up to rigid body motions for $\Pkm$.
\end{theorem}

Recall that $v^T \Big(\sum_{i=1}^m \lambda_i^* \nabla^2 f_i(p^*)\Big) v=v^T \nabla^2_p L(p^*,\lambda^*) v$.

This theorem is closely related to the standard theorem regarding second-order sufficient conditions in \cite{nocedal1999numerical}, Theorem 12.6, except it is slightly modified to deal with rigid-body motions. For this reason, we provide a proof in the Appendix \ref{app:proof-thm216}.

\section{The connection between rigidity and optimization}\label{sec:our-theorems}

In this section we show the connection between points $p^*$ which satisfy various necessary or sufficient conditions for the optimization problems $\Pkpm$ in \eqref{eqn:opt}, and rigidity properties of the associated framework $(p^*,E)$ (Section \ref{sec:reg}). 
We also consider a variation of $\Pkpm$ which allows designing certain components of the stress (Section \ref{sec:stress}).  

\subsection{The connection between solutions to \texorpdfstring{$\Pkpm$}{Pk plus or minus} and prestress stable frameworks}\label{sec:reg}

Our main results can be understood from the formal similarity between optimization conditions for $\Pkpm$, and conditions for various flavors of rigidity. 

Consider first the KKT conditions \eqref{eqn:KKT} for $\Pkpm$, evaluated at point $p$. These can be written in terms of the rigidity matrix $R(p)$ (see \eqref{R}) as
\begin{equation*}
	\lambda^T R(p)= 0.
\end{equation*}
Hence, we see immediately that when the KKT conditions hold, then the Lagrange multiplier $\lambda$ is a self-stress for framework $(p,E)$. Conversely, if $(p,E)$ has a self-stress $\omega$ with $\omega_k \neq 0$, then we obtain a Lagrange multiplier by scaling the $k$th coordinate to equal 1, as $\lambda = \omega/\omega_k$. We thus obtain a key result:

\begin{lemma}\label{lem:KKT}
	The point $p$ satisfies the KKT conditions \eqref{eqn:KKT} for optimization problem $\Pkpm$  \eqref{eqn:opt} with Lagrange multiplier $\lambda$, if and only if the framework $(p,E)$  has a self-stress $\omega$ with $\omega_k \neq 0$. 
	Furthermore, the Lagrange multipliers are in one-to-one correspondence with the self-stresses up to scaling: each self-stress $\omega$ gives rise to a Lagrange multiplier  $\lambda = \omega/\omega_k$, and each Lagrange multiplier $\lambda$ gives rise to a direction of self-stresses spanned by $\omega=\lambda$.
\end{lemma}

It will be useful later to characterize the number of self-stresses, which is offered in this lemma. 

\begin{lemma}\label{lem:localLICQ}
	Suppose $p^*$ is a local solution  to one of $\Pkpm$ and the LICQ condition holds. Then framework $(p^*,E)$ has a one-dimensional space of self-stresses. 
\end{lemma}

\begin{proof}
	If $p^*$ is a local solution to one of $\Pkpm$ and the LICQ condition holds, then $p^*$ satisfies the KKT conditions (Theorem \ref{thm:firstordernecessary}), and the Lagrange multiplier  $\lambda^*$ is unique (Remark \ref{rmk:unique-lagrange-multiplier}). Therefore there is a one-dimensional space of self-stresses  (Lemma \ref{lem:KKT}). 
\end{proof}

Now we may show that a solution to $\Pkpm$ is rigid. 

\begin{theorem}[An optimal edge gives a rigid framework]\label{thm:rigidity-optimal}
	Suppose $p^*$ is a strict local solution up to rigid body motions to either $\Pkp$ or $\Pkm$. Then the framework $(p^*,E)$ is locally rigid.
	Furthermore, if $E$ is isostatic or under-constrained, then $(p^*,E)$ is not first-order rigid. 
\end{theorem}

Note that the implication does not go the other way: if $(p,E)$ is rigid, this does not imply that $p$ is a solution, strict or otherwise, to one of $\Pkpm$. This is clear since there are many frameworks that are first-order rigid. 

\begin{proof}
	Suppose $p^*$ is a strict local solution up to rigid body motions to either $\Pkp$ or $\Pkm$. Then there is a neighborhood $U$ of $p$ such that any $q\in U$ with $q$ not congruent to $p$, has $f_k(q)\neq f_k(p)$. Therefore, any $q\in U$ with $f_k(q)= f_k(p)$ must be congruent to $p$.  Take a point $q'\in U$ with the same edge lengths as $p$. Then $f_k(q')=f_k(p)$, so we must have that $q'$ is congruent to $p$. Hence, $p$ is locally rigid. 
	
	Now we show that if $E$ is isostatic or under-constrained, then $(p^*,E)$ is not first-order rigid. For this we must show there is a self-stress (Remark \ref{rmk:fundamental-lin-alg}). If the LICQ condition holds, then Lemma \ref{lem:localLICQ} shows the existence of a self-stress. If the LICQ condition doesn't hold, then there are constants $\{c_i\}_{i\neq k}$, at least two of which are nonzero, such that $\sum_{i\neq k}c_i\nabla f_i(p^*) = 0$. Defining vector $\omega\in \R^m$ with $\omega_k = 0$ and $\omega_i = c_i$ ($i\neq k$), we see immediately that $\omega$ is a self-stress.  
\end{proof}

Here is a technical lemma that we'll need later. 

\begin{lemma}\label{lem:spaces}
	Suppose $p^*$ is a feasible point for $\Pkpm$, which satisfies the KKT conditions \eqref{eqn:KKT} for some Lagrange multiplier $\lambda^*$. Then the following spaces are equivalent:
	\begin{equation}
		\nullspace \: R(p^*) \cap \mathcal{T}(p^*)^\perp = {C_k(p^*)}.
	\end{equation}
\end{lemma}

Recall that $C_k(p^*)$ was defined in \eqref{eqn:cone-rigid}.

\begin{proof}
	The definition of $\nullspace \: R(p^*)\cap \mathcal{T}(p^*)^\perp $ is identical to that of ${C_k(p^*)}$ except it includes the additional condition $\nabla f_k(p^*)v = 0$. 
	Therefore, by definition, we have $\nullspace \: R(p^*) \cap \mathcal{T}(p^*)^\perp \subset {C_k(p^*)}$. 
	To show the inclusion goes the other way, note that any $v$ such that $v \perp \nabla f_i(p^*)$ for all $i\neq k$ must also satisfy $v \perp \nabla f_k(p^*)$, by the KKT conditions. 
	Therefore $\nullspace \: R(p^*) \cap \mathcal{T}(p^*)^\perp \supseteq {C_k(p^*)}$.
\end{proof}

How rigid is a solution $p^*$ to the optimization problem \eqref{eqn:opt}? This depends on whether it satisfies the second-order sufficient conditions, which will guarantee that it is prestress stable. 
It is convenient to rewrite the second-order sufficient conditions \eqref{eqn:secondsufficient}  using rigidity theory notation. Observe that
\begin{equation}\label{eqn:vLv}
	\frac{1}{2}v^T \nabla^2_pL(p^*,\lambda^*) v = \sum_{i=1}^m \lambda_i^* |v_{i,1}-v_{i,2}|^2 = \lambda^T R(v)v.
\end{equation}
Hence, the inequality in \eqref{eqn:secondsufficient} can be  written as $\lambda^T R(v)v>0$, which can be compared with the inequality $\omega^TR(v)v>0$ in condition \eqref{eqn:prestress-stable} for prestress stability. 

\begin{theorem}[Second-order sufficient conditions imply prestress stability]\label{thm:prestress-stability}
	Suppose $p^*$ is a strict local solution for one of $\Pkpm$ which is certified by Theorem \ref{thm:secondordersufficient}: the KKT conditions \eqref{eqn:KKT} hold for some Lagrange multiplier $\lambda^*$, and the second-order sufficient conditions \eqref{eqn:secondsufficient} hold with the appropriate direction of inequality. That is, for every $v \in C_k(p^*)$ and $v \neq 0$,
	\begin{equation}\label{eqn:2nd-condition-pm}
		v^T \Big(\sum_{i=1}^m \lambda_i^* \nabla^2 f_i(p^*)\Big) v > 0 \quad \text{for }P_k^+,\quad v^T \Big(\sum_{i=1}^m \lambda_i^* \nabla^2 f_i(p^*)\Big) v < 0 \quad \text{for }P_k^-.
	\end{equation}
	Then the framework $(p^*,E)$ is prestress stable. If additionally  $E$ is isostatic or under-constrained, then $(p^*,E)$ is not first-order rigid. 
\end{theorem}

\begin{proof}
	First consider problem $\Pkp$. Under the conditions of the theorem there is a  self-stress $\omega=\lambda^*$ for the framework $(p^*,E)$ (Lemma \ref{lem:KKT}). 
	By \eqref{eqn:vLv}, 
	\begin{equation}\label{eqn:wRvv}
		\omega^TR(v)v > 0 \qquad \forall \; v \in {C_k(p^*)}, \quad v\neq 0.
	\end{equation}
	To show that $(p^*,E)$ is prestress stable, i.e. it satisfies \eqref{eqn:prestress-stable}, we must show the above inequality holds for all $v\in \nullspace \: R(p^*) \cap \mathcal{T}(p^*)^\perp$, $v\neq 0$. 
	When the KKT conditions hold, then  $\nullspace \: R(p^*) \cap \mathcal{T}(p^*)^\perp={C_k(p^*)}$ (Lemma \ref{lem:spaces}).
	Hence, the second-order sufficient conditions are equivalent to the condition for prestress stability, with the self-stress $\omega$ showing prestress stability, equal to the  Lagrange multiplier $\lambda^*$ for the KKT conditions. 
	
	If instead the second-order sufficient conditions for problem $\Pkm$ hold, then construct the self-stress for prestress stability as $\omega = -\lambda^*$. The remaining argument is identical. 
	
	If $E$ is isostatic or under-constrained, then since there exists a self-stress (Lemma \ref{lem:KKT}), framework $(p^*,E)$ cannot be first-order rigid. 
\end{proof}

We might ask if the implication goes the other way -- if we have a prestress stable framework, is it the solution to an optimization problem of the form $\Pkpm$ for some $k$? The answer is yes, for each edge such that the self-stress guaranteeing prestress stability is nonzero on that edge.

\begin{theorem}[Prestress stability implies edge lengths are optimal]\label{thm:2nd-sufficient-cond}
	Suppose framework $(p^*,E)$ is prestress stable, and let $\omega$ be the self-stress that guarantees prestress stability in \eqref{eqn:prestress-stable}. Then, for each edge $E_k$ such that $\omega_k>0$, $p^*$ is a strict local solution up to rigid body motions for $\Pkp$, and for each edge $E_k$ such that $\omega_k<0$, $p^*$ is a strict local solution up to rigid body motions for $\Pkm$.
\end{theorem}

\begin{proof}
	Let $E_k$ be an edge such that $\omega_k\neq 0$. Then $p^*$ satisfies the KKT conditions with Lagrange multiplier $\lambda = \omega/\omega_k$ (Lemma \ref{lem:KKT}). 
	By the condition for prestress stability, we have that 
	\[
	\omega^TR(v)v >0\qquad  \text{for } v\in \nullspace \: R(p^*) \cap \mathcal{T}(p^*)^\perp, \quad v\neq 0.
	\]
	Lemma \ref{lem:spaces} shows that $\nullspace \: R(p^*) \cap \mathcal{T}(p^*)^\perp={C_k(p^*)}$. Hence, by \eqref{eqn:vLv}, if $\omega_k >0$, using $\lambda = \omega/\omega_k$, we have
	\[
	\lambda^TR(v)v >0,\qquad  \text{for } v\in {C_k(p^*)}, \quad v\neq 0.
	\]
	Therefore the second-order conditions \eqref{eqn:secondsufficient} are verified, so, by Theorem \ref{thm:secondordersufficient}, $p^*$ is a strict local solution up to rigid body motions for $\Pkp$. 
	
	If instead $\omega_k <0$, then we have 
	\[
	\lambda^TR(v)v <0,\qquad  \text{for } v\in {C_k(p^*)}, \quad v\neq 0.
	\]
	Hence,  by Theorem \ref{thm:secondordersufficient}, $p^*$ is a strict local solution up to rigid body motions for $\Pkm$. 
\end{proof}

Theorems \ref{thm:prestress-stability} and \ref{thm:2nd-sufficient-cond} give us a new way to characterize a prestress stable framework. 
Suppose we have a prestress stable framework $(p,E)$, with self-stress $\omega$ certifying prestress stability. Then \emph{every edge $E_k$ where $\omega_k\neq 0$ is individually at a local optimum}! That is, is we freeze the lengths of every edge except $E_k$, and ask to minimize the length of $E_k$ (if $\omega_k >0$) or maximize the length of $E_k$ (if $\omega_k <0$), then we can go no further -- edge $E_k$ is already optimized!

Conversely, suppose we solve an optimization problem of the form $\Pkp$, obtaining a strict local solution certified by the second-order sufficient conditions, and hence a prestress stable framework (Theorem \ref{thm:prestress-stability}). Then \emph{every other edge where the Lagrange multiplier is nonzero, is individually at a local optimum!} Therefore by optimizing the length of just one edge, we simultaneously optimize the lengths of all other edges. 

\subsection{Stress design}\label{sec:stress}

We now consider a variation of $\Pkpm$, which will allow us to set components of the stress. 
To motivate the variation, suppose we replace the objective function in $\Pkpm$ with $f_1(p) + f_2(p)$. Then, the KKT conditions for this modified problem would ask for a set of numbers $\lambda_3,\dots,\lambda_m$ such that 
\begin{equation*}
	\nabla f_1(p) + \nabla f_2(p) + \lambda_3\nabla f_3(p) + \ldots + \lambda_m \nabla f_m(p) = 0,
\end{equation*}
which can be written using rigidity notation as 
\begin{equation*}
	(1,1,\lambda_3,\dots,\lambda_m)^T \: R(p) = 0.
\end{equation*}
Defining the ``Lagrange multiplier'' to be $\lambda = (1,1,\lambda_3,\dots,\lambda_m)$, we see that the existence of a Lagrange multiplier $\lambda$ implies the existence of a self-stress $\omega=\lambda$. In this self-stress, the ratio of the stress on edges 1 and 2 is fixed at 1:1. 

In general, we can consider a linear combination of edge-length functions in our objective function, as $\alpha f_1(p) + \beta f_2(p)$, leading to a stress ratio of $\alpha:\beta$ on edges 1,2. We can also include more edges, thereby fixing the stress ratios on multiple edges. 

For the general problem, suppose $S\subset E$ is a set of edges whose stress ratios we wish to fix, and let $S^c$ be its complement in $E$. We wish to solve the following optimization problem. 

\begin{equation}\label{eqn:stress-design}
	\min \sum_{k\in S} \sigma_k f_k  \qquad \text{s.t.} \qquad f_j(p) - l_j^2 = 0, \qquad j \in S^c.
\end{equation}
That is, for each edge, we either fix its length, or fix the component of the stress on that edge.

We now sketch how our results extend to this case; the rigorous proofs are nearly identical to those in Section \ref{sec:reg}. First, the KKT conditions ask that 
\begin{equation*}
	\sum_{k\in S} \sigma_k \nabla f_k + \sum_{k\in S^c} \lambda_k \nabla f_k = 0. 
\end{equation*}
This can be written as $\lambda^T R(p)=0$ where we define $\lambda_k = \sigma_k$ if $k\in S$. Hence, the extension of Lemma \ref{lem:KKT} holds: a Lagrange multiplier $\lambda$ induces a self-stress $\omega=\lambda$, and conversely, given a self-stress $\omega$ with $\omega_k = \sigma_k$ for $k\in S$, the KKT conditions hold with $\lambda = \omega$. Therefore, a local solution $p^*$ to \eqref{eqn:stress-design} gives a framework $(p^*,E)$ with the self-stresses fixed (up to scaling) on the edges in $S$. 

Additionally, the second-order sufficient condition asks that
\begin{align*}
	v^T \Big(\sum_{i=1}^m \lambda_i \nabla^2 f_i(p^*)\Big) v > 0, \quad \forall v\in C_S(p^*) \quad v \neq 0,
\end{align*}
where the critical cone $C_S(p^*)$ in this case is
\begin{align*}
	C_S(p^*) = \{v \in \mathcal{T}(p^*)^\perp \:|\: \nabla f_j(p^*) v = 0, \quad j \in S^c\}.
\end{align*}
With the KKT and second-order sufficient condition, the extension of Theorem \ref{thm:prestress-stability} holds, i.e. the corresponding bar framework is prestress stable since $\nullspace \: R(p^*) \cap \mathcal{T}(p^*)^\perp \subset C_S(p^*)$. However, Theorem \ref{thm:2nd-sufficient-cond} extends only in the special case  where $\nullspace R(p^*) \cap \mathcal{T}(p^*)^\perp = C_S(p^*)$. \footnote{The equality holds when $\{\nabla f_i(p)\}_{i \in S}$ are in the subspace generated by $\mbox{span}\{\nabla f_j(p),{j \in S^c}\}$. This holds automatically when $|S|=1$, because the Lagrange multiplier ensures that $\nabla f_i$ is spanned by $\{\nabla f_j\}_{j\neq i}$ as in Lemma \ref{lem:spaces}. However when there is more than one free edge, the existence of a Lagrange multiplier does not automatically guarantee that each individual gradient is in this span.}

We remark that while the KKT conditions ensure the existence of one self-stress with the desired stress ratios on the edges, there could be other self-stresses whose stress ratios are uncontrolled. 

\begin{remark}\label{rmk:force-density}
	A common stress design method used in the engineering literature is the force density method \cite{schek1974force}. For a $d$-dimensional bar framework with $n$ vertices and $m$ edges, the force density method directly solves the linear system $w^T R(p) = 0$ for $p\in \mathbb{R}^{nd}$ with the given self-stress $w \in \mathbb{R}^m$. In fact, the force density method directly corresponds to the unconstrained optimization problem with objective function $\sum_{i=1}^m w_i |p_{i,1}-p_{i,2}|^2$ and no constraints on the edge lengths.
\end{remark}

\section{Designing prestress stable frameworks using optimization}\label{sec:examples}

Theorem \ref{thm:prestress-stability}  immediately suggests an algorithm for designing frameworks that are prestress stable but not first-order rigid: start with an isostatic or under-constrained framework $(p,E)$, and apply an optimization algorithm to minimize or maximize the length of some edge $E_k$, by solving one of $\Pkpm$ in \eqref{eqn:opt}. 

This approach is guaranteed to gives us two different prestress stable frameworks (one for the solution to each of $\Pkpm$), under weak conditions:  the length of $E_k$ is not constant on the feasible set, and, for $\Pkp$, this length is bounded away from 0 on the feasible set. 

In this section we show several examples of frameworks designed by this method\footnote{Animations corresponding to all figures in this paper are available at\newline\url{https://github.com/xuenanli/constrained-optimization-rigidity}.}. We note that optimization applied to most ``typical'' graphs gives vertices that are collinear or coplanar, which we find uninteresting. Therefore, to construct interesting examples, where the existence of self-stresses and $\mathcal{T}^{\perp}$-flexes is not obvious upon visual inspection, we do the following:
\begin{itemize}
	\item We start from a known prestress stable framework that we find interesting, and randomly perturb the positions of the vertices. This makes the framework first-order rigid and  the edges have random lengths. 
	
	\item We select one edge for our constrained optimization and numerically minimize or maximize its length, while constraining the lengths of all other edges. 
\end{itemize}

To eliminate the effect of rigid body motions, we use a pinning scheme which fixes the value of $D := d(d + 1)/2$ coordinates to zero. Details on how to choose the coordinates, and a brief explanation of why this leaves the rigidity properties invariant, are given in Appendix~\ref{app:pinning}.

We use a projected gradient descent method to numerically search for local solutions of 
\eqref{eqn:opt} with additional \edit{pinning constraints}
(the details of our numerical scheme are provided in Appendix~\ref{app:numerical-scheme}). A numerical solution $p^*$
is prestress stable if it satisfies the KKT condition and the second-order sufficient condition (adapted for pinning constraints in \eqref{eqn:2nd-test-pinning}).

\begin{figure}[!htb]
	\centering
    \tikzset{
        pinned flex/.style={-{Latex[length=1.5mm,width=1mm]},thick,
            shorten <=0.6mm,preaction={draw=white,line width=1.6pt,-}},
        pin mark/.style={black,line width=0.8pt,
            preaction={draw=white,line width=1.8pt}}
    }
	\subfloat[]{
		\begin{tikzpicture}[scale=0.7]
			\coordinate (A1) at (0, 0);
			\coordinate (B1) at (-5.0000e-01, 8.6603e-01);
			\coordinate (C1) at (6.2583e-01, 1.5160e+00);
			\coordinate (D1) at (2.6258e+00, 1.5160e+00);
			\coordinate (E1) at (3.4744e+00, 6.6750e-01);
			\coordinate (F1) at (1, 0);
			
			\draw[ultra thick] (A1) -- (B1) -- (C1) -- (D1) -- (E1) -- (F1) -- cycle;
			
			\draw[ultra thick] (B1) -- (E1);
			\draw[ultra thick] (C1) -- (F1);
			\draw[red, ultra thick,dashed] (A1) -- (D1);
			
			\draw[pin mark] ($(A1)+(-4pt,-4pt)$) -- ($(A1)+(4pt,4pt)$)
                ($(A1)+(-4pt,4pt)$) -- ($(A1)+(4pt,-4pt)$);
			\draw[pin mark] ($(F1)+(-4pt,-4pt)$) -- ($(F1)+(4pt,4pt)$)
                ($(F1)+(-4pt,4pt)$) -- ($(F1)+(4pt,-4pt)$);
			\fill[black] (B1) circle (3pt);
			\fill[black] (C1) circle (3pt);
			\fill[black] (D1) circle (3pt);
			\fill[black] (E1) circle (3pt);
		\end{tikzpicture}
	}
	\hfil
	\subfloat[]{
		\begin{tikzpicture}[scale=0.7]
			\coordinate (A2) at (0,0);
			\coordinate (B2) at (-6.1425e-01, 7.8911e-01);
			\coordinate (C2) at (4.8962e-01, 1.4758e+00);
			\coordinate (D2) at (2.4634e+00, 1.7987e+00);
			\coordinate (E2) at (3.3594e+00, 1.0005e+00);
			\coordinate (F2) at (1, 0);
			
			\draw[ultra thick,blue] (A2) -- (B2); 
			\draw[ultra thick,blue] (B2) -- (C2); 
			\draw[ultra thick,blue] (C2) -- (D2); 
			\draw[ultra thick,blue] (D2) -- (E2); 
			\draw[ultra thick,blue] (E2) -- (F2); 
			\draw[ultra thick,blue] (A2) -- (F2); 
			\draw[ultra thick,red] (B2) -- (E2); 
			\draw[ultra thick,red] (C2) -- (F2); 
			\draw[ultra thick, red] (A2) -- (D2); 
			
			\draw[pin mark] ($(A2)+(-4pt,-4pt)$) -- ($(A2)+(4pt,4pt)$)
                ($(A2)+(-4pt,4pt)$) -- ($(A2)+(4pt,-4pt)$);
			\draw[pin mark] ($(F2)+(-4pt,-4pt)$) -- ($(F2)+(4pt,4pt)$)
                ($(F2)+(-4pt,4pt)$) -- ($(F2)+(4pt,-4pt)$);
			\fill[black] (B2) circle (3pt);
			\fill[black] (C2) circle (3pt);
			\fill[black] (D2) circle (3pt);
			\fill[black] (E2) circle (3pt);
			
            \def\flexscale{-1.5}
            \foreach \vertex/\vx/\vy in {
                B2/0.220672338748/0.171773187717,
                C2/0.269547096380/0.093218865074,
                D2/0.367032129902/-0.502668396638,
                E2/0.262776876285/-0.619687428348}{
                \draw[pinned flex] (\vertex) -- ++({\flexscale*\vx},{\flexscale*\vy});
            }
		\end{tikzpicture}
	}
	\hfil
	\subfloat[]{
		\begin{tikzpicture}[scale=0.7]
			\coordinate (A1) at (0, 0);
			\coordinate (B1) at (6.5000e-01, 1.1258e+00);
			\coordinate (C1) at (1.6500e+00, 1.1303e+00);
			
			\coordinate (D1) at (3.2000e+00, 0);
			\coordinate (E1) at (3.4568e+00, 2.1850e+00);
			\coordinate (F1) at (2.3762e+00, 8.7259e-01);
			
			\draw[ultra thick] (A1) -- (B1);
			\draw[ultra thick] (B1) -- (C1);
			\draw[ultra thick] (A1) -- (C1);
			
			\draw[ultra thick] (D1) -- (E1);
			\draw[ultra thick] (E1) -- (F1);
			\draw[ultra thick] (D1) -- (F1);
			
			\draw[ultra thick] (B1) -- (E1);
			\draw[ultra thick] (A1) -- (D1);
			\draw[red, ultra thick,dashed] (C1) -- (F1);
			
			\draw[pin mark] ($(A1)+(-4pt,-4pt)$) -- ($(A1)+(4pt,4pt)$)
                ($(A1)+(-4pt,4pt)$) -- ($(A1)+(4pt,-4pt)$);
			\draw[pin mark] ($(D1)+(-4pt,-4pt)$) -- ($(D1)+(4pt,4pt)$)
                ($(D1)+(-4pt,4pt)$) -- ($(D1)+(4pt,-4pt)$);
			
			\fill[black] (B1) circle (3pt);
			\fill[black] (C1) circle (3pt);
			\fill[black] (E1) circle (3pt);
			\fill[black] (F1) circle (3pt);
		\end{tikzpicture}
	}
	\hfil
	\subfloat[]{
		\begin{tikzpicture}[scale=0.7]
			\coordinate (A2) at (0, 0);
			\coordinate (B2) at (8.5812e-01, 9.7654e-01);
			\coordinate (C2) at (1.8396e+00, 7.8484e-01);
			
			\coordinate (D2) at (3.2000e+00, 0);
			\coordinate (E2) at (3.6146e+00, 2.1606e+00);
			\coordinate (F2) at (2.4417e+00, 9.3002e-01);
			
			\draw[ultra thick,blue] (A2) -- (B2); 
			\draw[ultra thick, red] (B2) -- (C2); 
			\draw[ultra thick,red] (A2) -- (C2); 
			
			\draw[ultra thick,blue] (D2) -- (E2); 
			\draw[ultra thick,red] (E2) -- (F2); 
			\draw[ultra thick,red] (D2) -- (F2); 
			
			\draw[ultra thick,blue] (B2) -- (E2); 
			\draw[ultra thick, blue] (A2) -- (D2); 
			\draw[red, ultra thick] (C2) -- (F2); 
			
			\draw[pin mark] ($(A2)+(-4pt,-4pt)$) -- ($(A2)+(4pt,4pt)$)
                ($(A2)+(-4pt,4pt)$) -- ($(A2)+(4pt,-4pt)$);
			\draw[pin mark] ($(D2)+(-4pt,-4pt)$) -- ($(D2)+(4pt,4pt)$)
                ($(D2)+(-4pt,4pt)$) -- ($(D2)+(4pt,-4pt)$);
			
			\fill[black] (B2) circle (3pt);
			\fill[black] (C2) circle (3pt);
			\fill[black] (E2) circle (3pt);
			\fill[black] (F2) circle (3pt);
			
            \def\flexscale{1.5}
            \foreach \vertex/\vx/\vy in {
                B2/0.389664856312/-0.342414642833,
                C2/0.313170476365/-0.734039007110,
                E2/0.264370122695/-0.050724730806,
                F2/0.113795831810/0.092792507654}{
                \draw[pinned flex] (\vertex) -- ++({\flexscale*\vx},{\flexscale*\vy});
            }
		\end{tikzpicture}
	}
	\caption{Isostatic, prestress stable frameworks designed by our optimization method. (a),(c) show first-order rigid frameworks used as initial conditions for optimization.  The red dashed edge is the chosen free edge. (b),(d) show the prestress stable frameworks found by maximizing (b) and minimizing (d) the free edges. The edge color indicates the sign of the self-stress: red means the edge is stretched ($\omega_k >0$) and blue means the edge is compressed ($\omega_k<0$). \edit{Crosses mark the pinned vertices, and the arrows show infinitesimal flexes satisfying the pinning constraints.}}
	\label{fig:2D-prestress-example}
\end{figure}
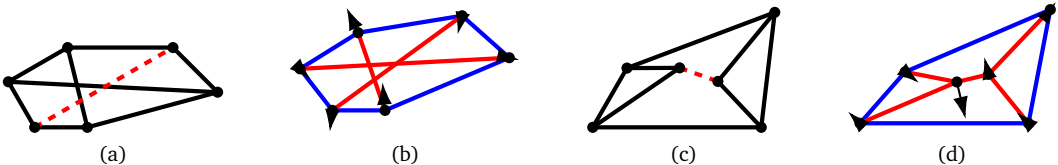

\paragraph{Isostatic examples}
Figure \ref{fig:2D-prestress-example} shows
two examples of two-dimensional isostatic frameworks that are prestress stable but not first-order rigid, which we henceforth refer to simply as \emph{prestress stable}. The first example features a perturbed hexagon with connected diagonals. By maximizing the length of the dashed edge in \cref{fig:2D-prestress-example}(a), we obtain the prestress stable framework in \cref{fig:2D-prestress-example}(b), which has \edit{one independent pinned infinitesimal flex} and one self-stress. The second example consists of two rigid triangles connected by three bars. Minimizing the length of the dashed edge in \cref{fig:2D-prestress-example}(c) results in another prestress stable framework with \edit{one independent pinned infinitesimal flex} and one self-stress. Notably, the configurations in \cref{fig:2D-prestress-example}(b) and (d) lack apparent symmetries, such as pairs of parallel lines or collinear points, and all the edges have different lengths.

\begin{figure}[!htb]
	\centering
	\captionsetup[subfloat]{farskip=-2.5pt,captionskip=-9.5pt}
	\subfloat[]{
		\includegraphics[width=0.2352\linewidth]{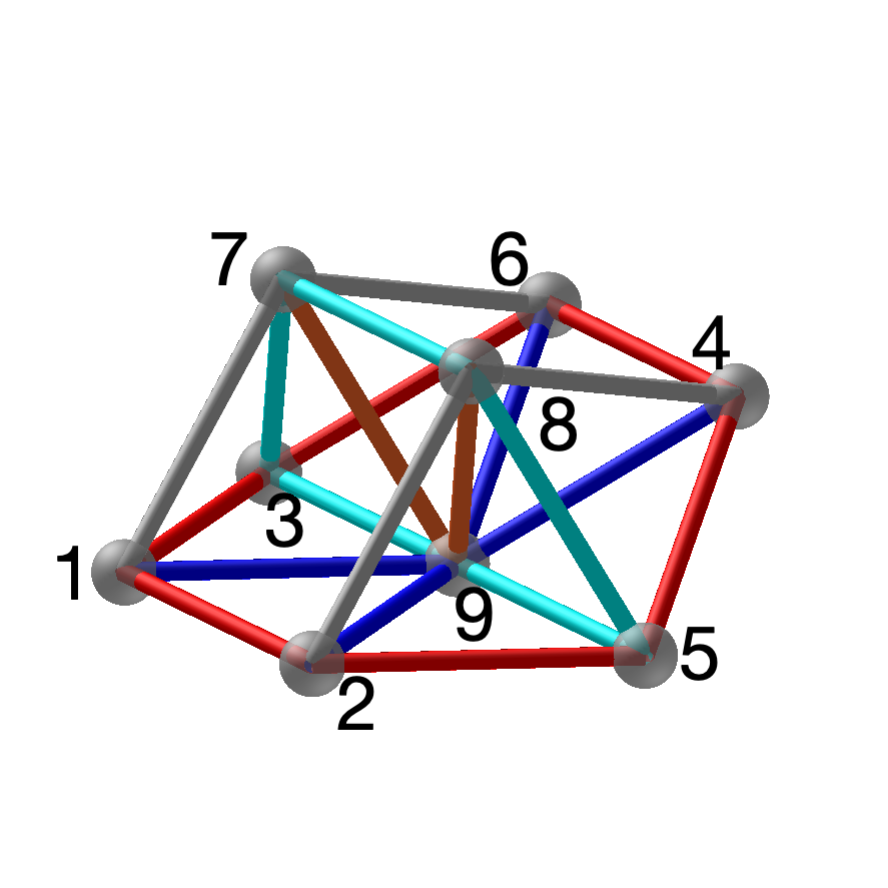}
	}
	\hfil
	\subfloat[]{
		\includegraphics[width=0.2352\linewidth]{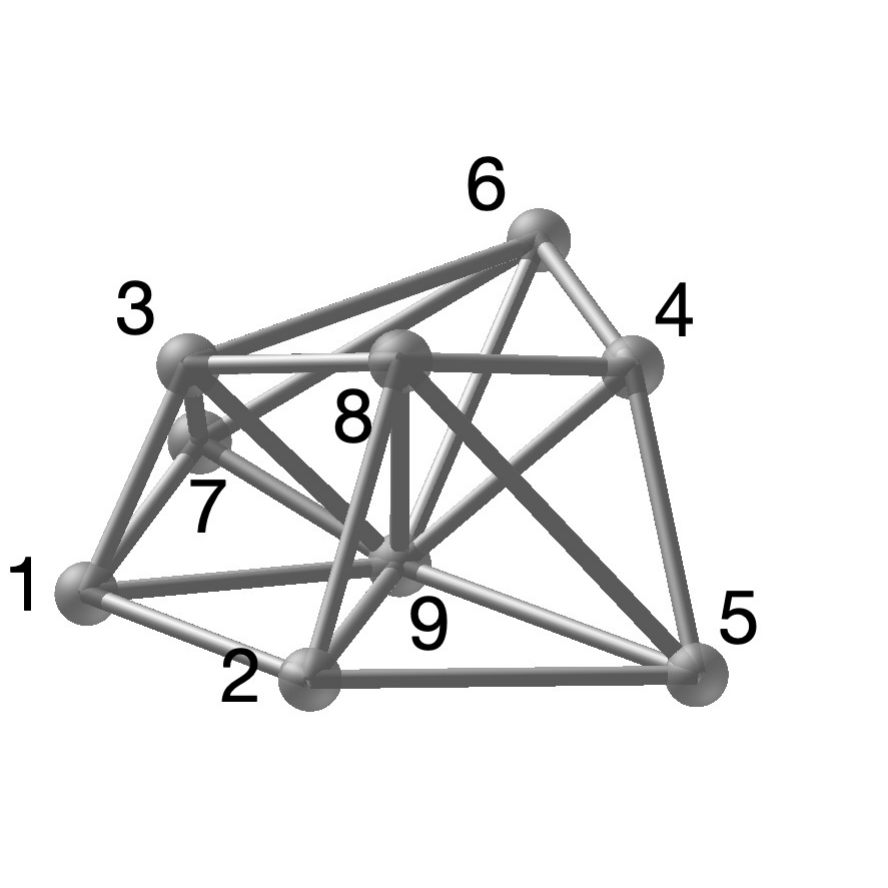}
	}
	\hfil
	{\captionsetup[subfloat]{captionskip=-7.5pt}
	\subfloat[]{
		\includegraphics[width=0.1848\linewidth]{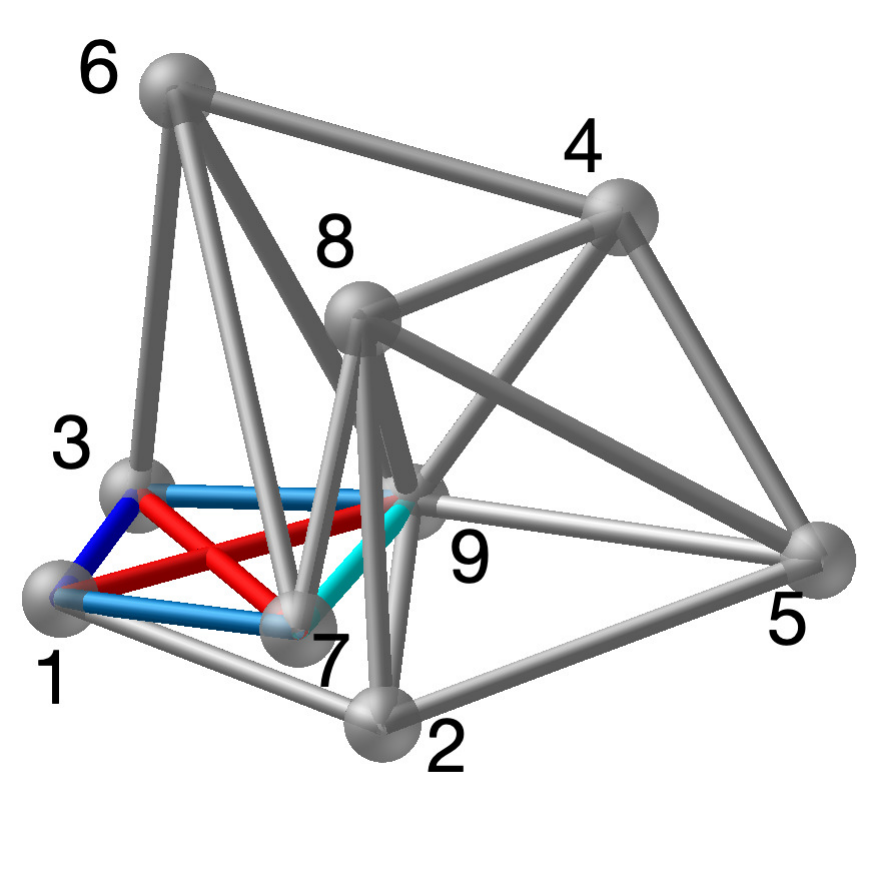}
	}}\\
	\subfloat[]{
		\includegraphics[width=0.2352\linewidth]{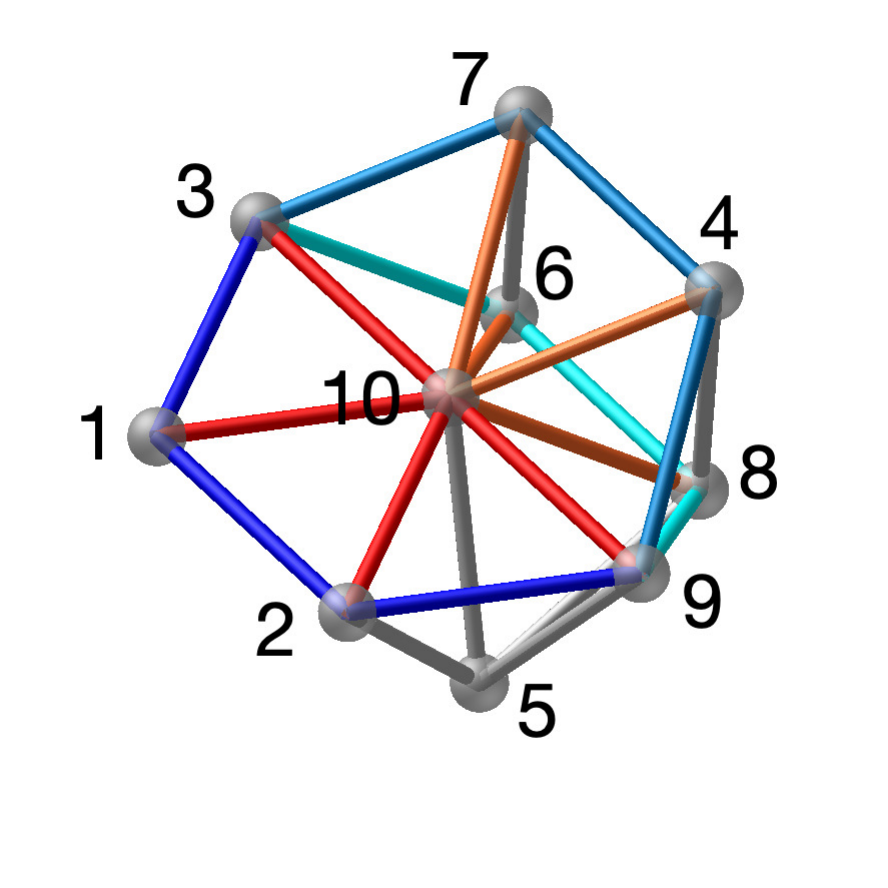}
	}
	\hfil
	\subfloat[]{
		\includegraphics[width=0.252\linewidth]{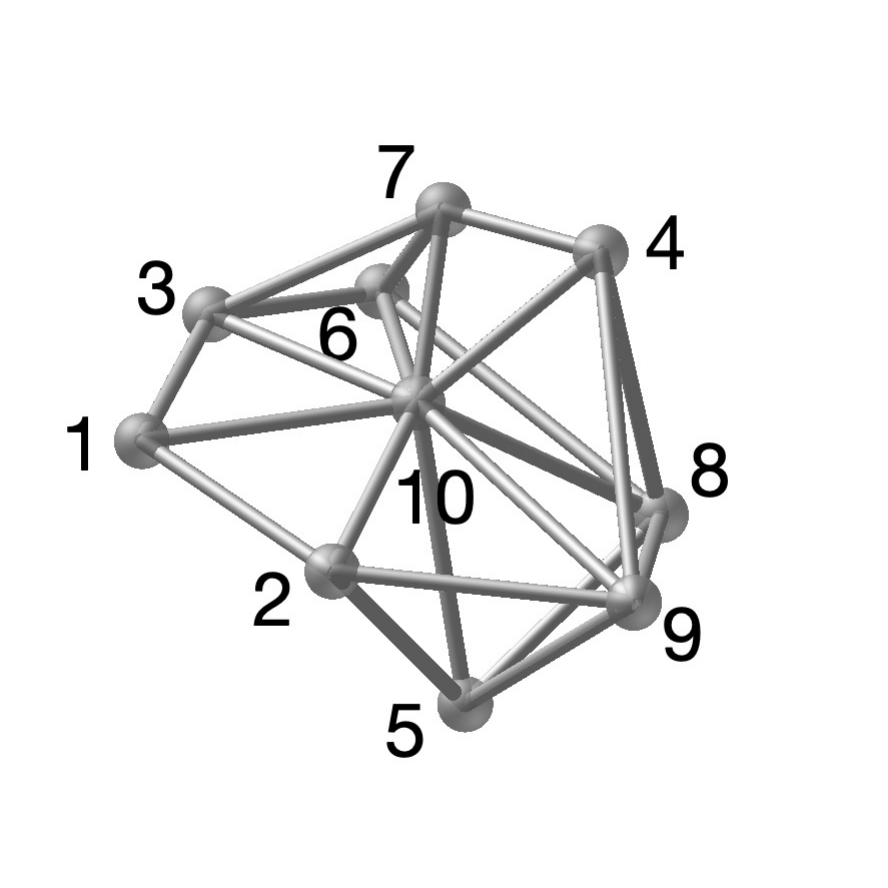}
	}
	\hfil
	\subfloat[]{
		\includegraphics[width=0.252\linewidth]{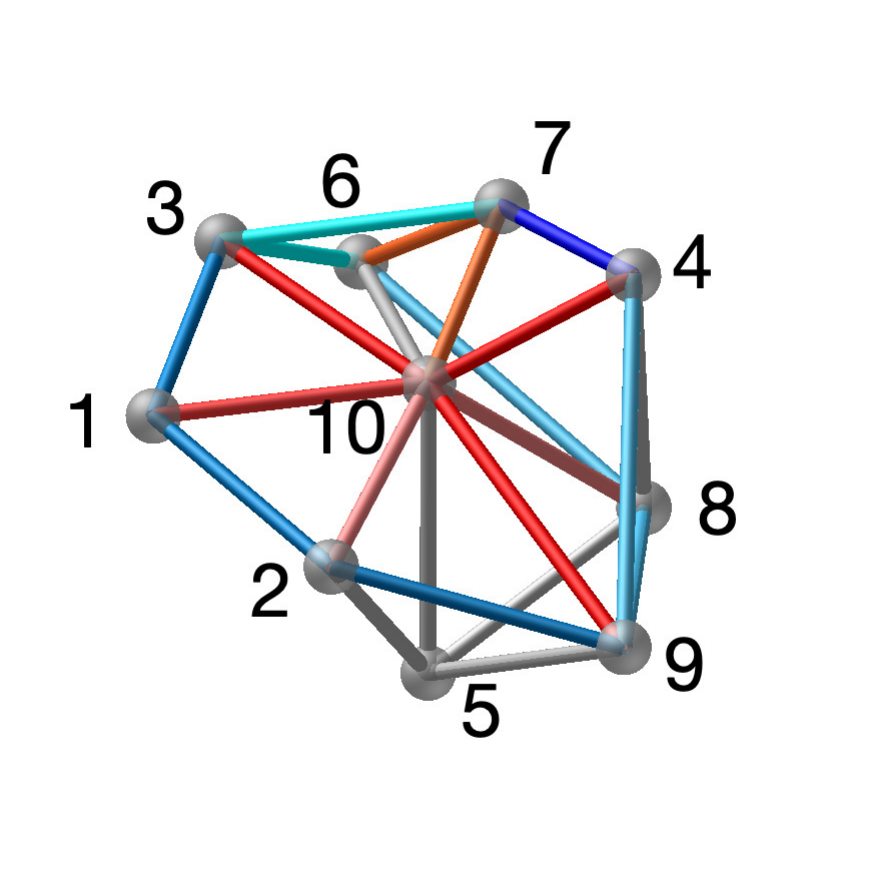}
	}
	\caption{3d prestress stable frameworks designed by our optimization method. (a) and (c) come from the sphere packing data; (b) and (e) are a perturbed version of (a) and (c); (c) and (f) are the prestress stable bar frameworks obtained by our constrained optimization scheme. Colors are proportional to the self-stress: red means under compression and blue means under tension. The color intensity represents the magnitude of the stress, with darker colors indicating higher stress values and lighter colors indicating lower stress values.}
	\label{fig:3d-opt-length}
\end{figure}

\paragraph{3d examples}
We apply the same scheme to some three-dimensional frameworks. Specifically, we selected two frameworks that are not first-order rigid from a dataset of frameworks obtained from packings of identical spheres \cite{holmes2016enumerating}. The framework in \cref{fig:3d-opt-length}(a) is isostatic and the framework in \cref{fig:3d-opt-length}(d) is under-constrained\footnote{Figure \ref{fig:3d-opt-length}(a) is the first case in $N=9$ and \cref{fig:3d-opt-length}(d) is the 6th case in $N=10$, which can be found in \url{https://personal.math.ubc.ca/~holmescerfon/packings.html}.}. We perturb the vertices randomly to obtain the frameworks in \cref{fig:3d-opt-length}(b) and \cref{fig:3d-opt-length}(e). By maximizing the distance between vertices 3 and 7 in \cref{fig:3d-opt-length}(b) and minimizing the distance between vertices 1 and 3 in \cref{fig:3d-opt-length}(e), we obtain two prestress stable frameworks with one self-stress in \cref{fig:3d-opt-length}(c) and \cref{fig:3d-opt-length}(f). Unlike the 2D examples, the 3D prestress stable frameworks obtained from our constrained optimization approach exhibit special geometry -- there are four vertices coplanar (vertices 1,3,7,9 for \cref{fig:3d-opt-length}(c) and vertices 1,3,7,10 for \cref{fig:3d-opt-length}(f)). 
Notably, these frameworks also have all edges different lengths. 

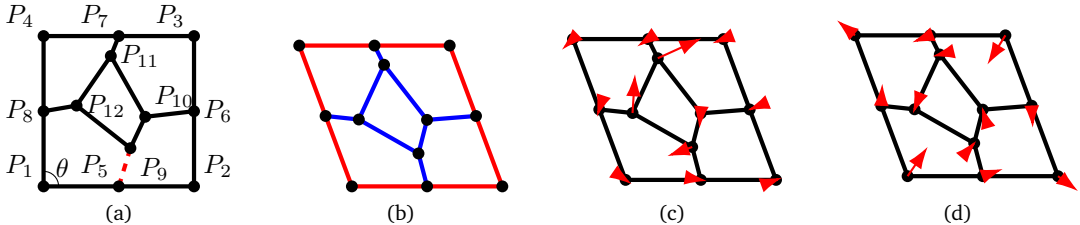
\begin{figure}[!htb]
	\centering
    \newcommand{\drawflexarrow}[3]{%
        \coordinate (direction) at ($(#1) + (#2,#3)$);
        \pgfmathsetmacro{\tipscale}{min(1,3.6*veclen(#2,#3))}
        \draw[red,line width=0.45pt,-{Latex[length=1.5mm,width=1mm,scale=\tipscale]},
            preaction={draw=white,line width=1.2pt,-}]
            ($(#1)!0.16cm!(direction)$) -- ++({2*#2},{2*#3});
    }
    \newcommand{\markflexpin}[1]{%
        \draw[black,line width=0.8pt,preaction={draw=white,line width=1.8pt}]
            ($(#1)+(-6pt,-6pt)$) -- ($(#1)+(6pt,6pt)$)
            ($(#1)+(-6pt,6pt)$) -- ($(#1)+(6pt,-6pt)$);
    }
    \subfloat[]{
		\begin{tikzpicture}[scale=0.4]
			\coordinate (A1) at (0, 0);
			\coordinate (A2) at (4,0);
			\coordinate (A3) at (4,4);
			\coordinate (A4) at (0,4);
			
			\coordinate (L1) at (2,0);
			\coordinate (L2) at (4,2);
			\coordinate (L3) at (2,4);
			\coordinate (L4) at (0,2);
			
			\coordinate (B1) at (2.3088e+00, 1.0172e+00);
			\coordinate (B2) at (2.7017e+00, 1.8515e+00);
			\coordinate (B3) at (1.7877e+00, 3.4701e+00);
			\coordinate (B4) at (8.7688e-01, 2.1496e+00);
			
			\draw[ultra thick,black] (A1) -- (A2); 
			\draw[ultra thick,black] (A2) -- (A3); 
			\draw[ultra thick,black] (A3) -- (A4); 
			\draw[ultra thick,black] (A4) -- (A1); 
			
			\draw[ultra thick,black] (B1) -- (B2); 
			\draw[ultra thick,black] (B2) -- (B3);
			\draw[ultra thick,black] (B3) -- (B4);
			\draw[ultra thick,black] (B4) -- (B1); 
			
			\draw[ultra thick,red,dashed] (L1) -- (B1); 
			\draw[ultra thick,black] (L2) -- (B2);
			\draw[ultra thick,black] (L3) -- (B3);
			\draw[ultra thick,black] (L4) -- (B4); 
			
			\draw (0,0) -- (0.4,0) arc[start angle=0, end angle=90, radius=0.4];
			\node at (0.5, 0.5) {$\theta$};
			
			\fill[black] (A1) circle (4.5pt);
			\node[above left] at (A1) {$P_1$};
			\fill[black] (A2) circle (4.5pt);
			\node[above right] at (A2) {$P_2$};
			
			\fill[black] (A3) circle (4.5pt);
			\node[above left] at (A3) {$P_3$};
			\fill[black] (A4) circle (4.5pt);
			\node[above left] at (A4) {$P_4$};
			
			\fill[black] (L1) circle (4.5pt);
			\node[above left] at (L1) {$P_5$};
			\fill[black] (L2) circle (4.5pt);
			\node[right] at (L2) {$P_6$};
			\fill[black] (L3) circle (4.5pt);
			\node[above left] at (L3) {$P_7$};
			\fill[black] (L4) circle (4.5pt);
			\node[left] at (L4) {$P_8$};
			
			\fill[black] (B1) circle (4.5pt);
			\node[below right] at (B1) {$P_9$};
			\fill[black] (B2) circle (4.5pt);
			\node[above right] at (B2) {$P_{10}$};
			\fill[black] (B3) circle (4.5pt);
			\node[right] at (B3) {$P_{11}$};
			\fill[black] (B4) circle (4.5pt);
			\node[right] at (B4) {$P_{12}$};
		\end{tikzpicture}
	}
	\hfil
	\subfloat[]{
		\begin{tikzpicture}[scale=0.4]
			\coordinate (A1) at (0, 0);
			\coordinate (A2) at (4,0);
			\coordinate (A3) at (2.6004e+00,3.7471e+00);
			\coordinate (A4) at (-1.3996e+00,3.7471e+00);
			
			\coordinate (L1) at (2,0);
			\coordinate (L2) at (3.3002e+00,1.8736e+00);
			\coordinate (L3) at (6.0036e-01,3.7471e+00);
			\coordinate (L4) at (-6.9982e-01,1.8736e+00);
			
			\coordinate (B1) at (1.7710e+00, 8.7585e-01);
			\coordinate (B2) at (1.9975e+00, 1.7698e+00);
			\coordinate (B3) at (8.5574e-01, 3.2366e+00);
			\coordinate (B4) at (1.8475e-01, 1.7795e+00);
			
			\draw[ultra thick,red] (A1) -- (A2); 
			\draw[ultra thick,red] (A2) -- (A3); 
			\draw[ultra thick,red] (A3) -- (A4); 
			\draw[ultra thick,red] (A4) -- (A1); 
			
			\draw[ultra thick,blue] (B1) -- (B2); 
			\draw[ultra thick,blue] (B2) -- (B3);
			\draw[ultra thick,blue] (B3) -- (B4);
			\draw[ultra thick,blue] (B4) -- (B1); 
			
			\draw[ultra thick,blue] (L1) -- (B1); 
			\draw[ultra thick,blue] (L2) -- (B2);
			\draw[ultra thick,blue] (L3) -- (B3);
			\draw[ultra thick,blue] (L4) -- (B4); 
			
			\fill[black] (A1) circle (4.5pt);
			\fill[black] (A2) circle (4.5pt);
			
			\fill[black] (A3) circle (4.5pt);
			\fill[black] (A4) circle (4.5pt);
			
			\fill[black] (L1) circle (4.5pt);
			\fill[black] (L2) circle (4.5pt);
			\fill[black] (L3) circle (4.5pt);
			\fill[black] (L4) circle (4.5pt);
			
			\fill[black] (B1) circle (4.5pt);
			\fill[black] (B2) circle (4.5pt);
			\fill[black] (B3) circle (4.5pt);
			\fill[black] (B4) circle (4.5pt);
		\end{tikzpicture}
	}
	\hfil
	\subfloat[]{
		\begin{tikzpicture}[scale=0.4]
			\coordinate (A1) at (0, 0);
			\coordinate (A2) at (4,0);
			\coordinate (A3) at (2.6004e+00,3.7471e+00);
			\coordinate (A4) at (-1.3996e+00,3.7471e+00);
			
			\coordinate (L1) at (2,0);
			\coordinate (L2) at (3.3002e+00,1.8736e+00);
			\coordinate (L3) at (6.0036e-01,3.7471e+00);
			\coordinate (L4) at (-6.9982e-01,1.8736e+00);
			
			\coordinate (B1) at (1.7710e+00, 8.7585e-01);
			\coordinate (B2) at (1.9975e+00, 1.7698e+00);
			\coordinate (B3) at (8.5574e-01, 3.2366e+00);
			\coordinate (B4) at (1.8475e-01, 1.7795e+00);
			
			\draw[ultra thick,black] (A1) -- (A2); 
			\draw[ultra thick,black] (A2) -- (A3); 
			\draw[ultra thick,black] (A3) -- (A4); 
			\draw[ultra thick,black] (A4) -- (A1); 
			
			\draw[ultra thick,black] (B1) -- (B2); 
			\draw[ultra thick,black] (B2) -- (B3);
			\draw[ultra thick,black] (B3) -- (B4);
			\draw[ultra thick,black] (B4) -- (B1); 
			
			\draw[ultra thick,black] (L1) -- (B1); 
			\draw[ultra thick,black] (L2) -- (B2);
			\draw[ultra thick,black] (L3) -- (B3);
			\draw[ultra thick,black] (L4) -- (B4); 
			
			\markflexpin{A1}
			\markflexpin{A2}
			
			\fill[black] (A3) circle (4.5pt);
			\fill[black] (A4) circle (4.5pt);
			
			\fill[black] (L1) circle (4.5pt);
			\fill[black] (L2) circle (4.5pt);
			\fill[black] (L3) circle (4.5pt);
			\fill[black] (L4) circle (4.5pt);
			
			\fill[black] (B1) circle (4.5pt);
			\fill[black] (B2) circle (4.5pt);
			\fill[black] (B3) circle (4.5pt);
			\fill[black] (B4) circle (4.5pt);
			
            \foreach \vertex/\vx/\vy in {
                A3/-0.32067420/-0.11977941,
                A4/-0.32067420/-0.11977941,
                L2/-0.16033710/-0.05988971,
                L3/-0.32067420/-0.11977941,
                L4/-0.16033710/-0.05988971,
                B1/-0.41128603/-0.10752701,
                B2/-0.15129626/-0.17339615,
                B3/0.34590041/0.21363854,
                B4/-0.10899028/0.42312042}{
                \drawflexarrow{\vertex}{\vx}{\vy}
            }
        \end{tikzpicture}
	}
	\hfil
	\subfloat[]{
		\begin{tikzpicture}[scale=0.4]
			\coordinate (A1) at (0, 0);
			\coordinate (A2) at (4,0);
			\coordinate (A3) at (2.6004e+00,3.7471e+00);
			\coordinate (A4) at (-1.3996e+00,3.7471e+00);
			
			\coordinate (L1) at (2,0);
			\coordinate (L2) at (3.3002e+00,1.8736e+00);
			\coordinate (L3) at (6.0036e-01,3.7471e+00);
			\coordinate (L4) at (-6.9982e-01,1.8736e+00);
			
			\coordinate (B1) at (1.7710e+00, 8.7585e-01);
			\coordinate (B2) at (1.9975e+00, 1.7698e+00);
			\coordinate (B3) at (8.5574e-01, 3.2366e+00);
			\coordinate (B4) at (1.8475e-01, 1.7795e+00);
			
			\draw[ultra thick,black] (A1) -- (A2); 
			\draw[ultra thick,black] (A2) -- (A3); 
			\draw[ultra thick,black] (A3) -- (A4); 
			\draw[ultra thick,black] (A4) -- (A1); 
			
			\draw[ultra thick,black] (B1) -- (B2); 
			\draw[ultra thick,black] (B2) -- (B3);
			\draw[ultra thick,black] (B3) -- (B4);
			\draw[ultra thick,black] (B4) -- (B1); 
			
			\draw[ultra thick,black] (L1) -- (B1); 
			\draw[ultra thick,black] (L2) -- (B2);
			\draw[ultra thick,black] (L3) -- (B3);
			\draw[ultra thick,black] (L4) -- (B4); 
			
			\markflexpin{A1}
			\markflexpin{A2}
			
			\fill[black] (A3) circle (4.5pt);
			\fill[black] (A4) circle (4.5pt);
			
			\fill[black] (L1) circle (4.5pt);
			\fill[black] (L2) circle (4.5pt);
			\fill[black] (L3) circle (4.5pt);
			\fill[black] (L4) circle (4.5pt);
			
			\fill[black] (B1) circle (4.5pt);
			\fill[black] (B2) circle (4.5pt);
			\fill[black] (B3) circle (4.5pt);
			\fill[black] (B4) circle (4.5pt);
			
            \foreach \vertex/\vx/\vy in {
                A3/-0.42544634/-0.15891177,
                A4/-0.42544634/-0.15891177,
                L2/-0.21272317/-0.07945588,
                L3/-0.42544634/-0.15891177,
                L4/-0.21272317/-0.07945588,
                B1/-0.14744817/-0.03854970,
                B2/-0.21739241/-0.02082795,
                B3/-0.33966377/-0.11600601,
                B4/-0.22234966/-0.17002957}{
                \drawflexarrow{\vertex}{\vx}{\vy}
            }
        \end{tikzpicture}
	}
	\caption{Designing a prestress stable framework with midpoint constraints. (a) the  initial condition, which is flexible; the dashed edge is the free edge whose length will be minimized; (b) the optimal state found by the constrained optimization approach. \edit{In (c-d), crosses mark the pinned vertices $P_1$ and $P_2$, and $P_5$ remains a midpoint of $P_1$ and $P_2$. The two arrow fields form a basis for the space of pinned infinitesimal flexes of (b).} \edit{The colors in (b) are the same as in \cref{fig:2D-prestress-example}.}}
	\label{fig:2D-square}
\end{figure}

\paragraph{Examples with additional linear constraints}
In \cref{fig:2D-square}, we explore a 2D example with a different kind of constraint, namely, we ask that four vertices remain at the midpoints of the four outer edges. 
Specifically, we examine a structure with an outer square and an inner quadrilateral in \cref{fig:2D-square}(a). Then we add four vertices ($P_5, P_6,P_7,P_8$) on the outer edges, connected with bars to the inner vertices, and we impose the constraint that these additional vertices remain at the midpoints of the edges they lie on (i.e. $P_5$ is halfway between $P_1,P_2$, etc). The whole system has 12 vertices (24 degrees of freedom), 12 edges, corresponding to (quadratic) distance constraints, and 8 (linear) midpoint constraints (2 for each vertex). This example is also explored in \cite{roback2025tuning}, where prestress stability arises from a special geometry (the inner quadrilateral is a square). In contrast, we use constrained optimization to obtain prestress stable structures without such special geometry.

The additional midpoint constraints are not distance constraints, so they don't fit the setup we have described in Section \ref{sec:prelim}. However, one can see by examining this section that the definitions and results are unchanged for any constraints $\{f_i\}$ that are invariant to translations and rotations. By defining the corresponding rigidity, i.e. local rigidity, first-order rigidity and prestress stability, one can achieve similar rigidity results shown in Section \ref{sec:reg}.

The initial framework in \cref{fig:2D-square} (a) is flexible with a mechanism that can be parametrized by an angle $\theta$ shown in \cref{fig:2D-prestress-example}(a). By minimizing the free edge (dashed in \cref{fig:2D-square}(a)), we obtain a prestress stable  framework \cref{fig:2D-square}(b). This framework could tile the plane, thus creating a lattice structure with a space of self-stresses. 

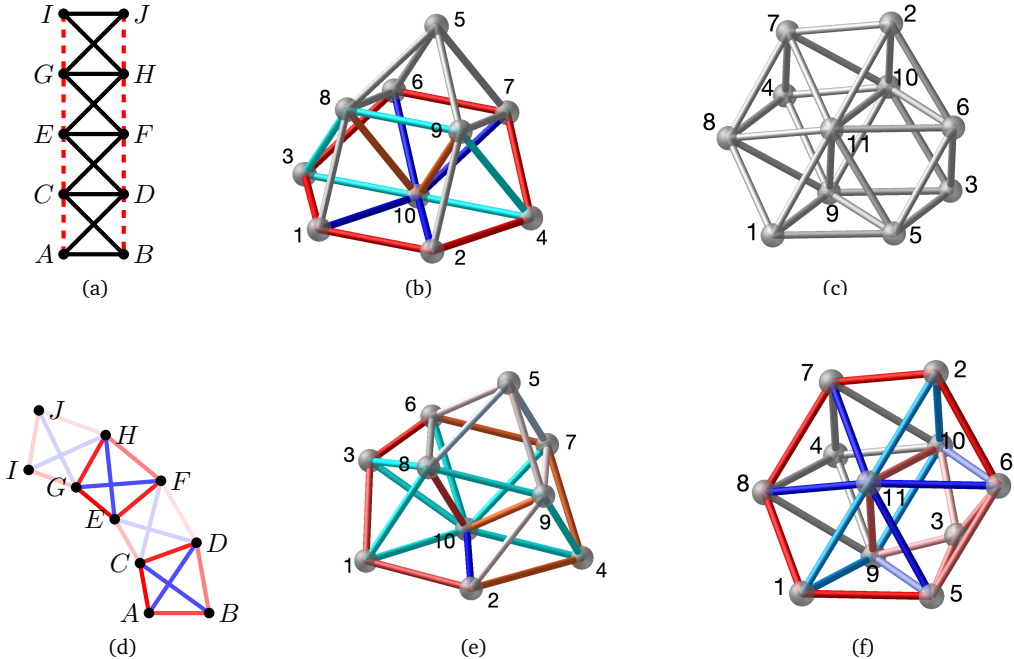
\begin{figure}[!htb]
	\centering
	\captionsetup[subfloat]{farskip=-3.5pt,captionskip=2pt}
	\subfloat[]{
		\begin{tikzpicture}[scale=0.5]
			\coordinate (A) at (0, 0);
			\coordinate (B) at (1, 0);
			\coordinate (C) at (0, 1);
			\coordinate (D) at (1, 1);
			\coordinate (E) at (0, 2);
			\coordinate (F) at (1, 2);
			\coordinate (G) at (0, 3);
			\coordinate (H) at (1, 3);
			\coordinate (I) at (0, 4);
			\coordinate (J) at (1, 4);
			
			\draw[ultra thick] (A) -- (B); 
			\draw[ultra thick, red, dashed] (A) -- (C); 
			\draw[ultra thick] (C) -- (D); 
			\draw[ultra thick, red, dashed] (D) -- (B); 
			\draw[ultra thick] (C) -- (B); 
			\draw[ultra thick] (A) -- (D); 
			
			\draw[ultra thick, red, dashed] (E) -- (C); 
			\draw[ultra thick] (E) -- (F); 
			\draw[ultra thick, red, dashed] (D) -- (F); 
			\draw[ultra thick] (C) -- (F); 
			\draw[ultra thick] (E) -- (D); 
			
			\draw[ultra thick, red, dashed] (E) -- (G); 
			\draw[ultra thick] (G) -- (H); 
			\draw[ultra thick, red, dashed] (F) -- (H); 
			\draw[ultra thick] (E) -- (H); 
			\draw[ultra thick] (G) -- (F); 
			
			\draw[ultra thick, red, dashed] (I) -- (G); 
			\draw[ultra thick] (I) -- (H); 
			\draw[ultra thick, red, dashed] (J) -- (H); 
			\draw[ultra thick] (G) -- (J); 
			\draw[ultra thick] (J) -- (I); 
			
			\foreach \point in {A, C, E, G, I} {
				\fill[black] (\point) circle (2.5pt);
				\node[left] at (\point) {$\point$};
			}
			
			\foreach \point in {B, D, F, H, J} {
				\fill[black] (\point) circle (2.5pt);
				\node[right] at (\point) {$\point$};
			}
		\end{tikzpicture}
	}
	\hfil
	\subfloat[]{
		\includegraphics[width=0.2\linewidth]{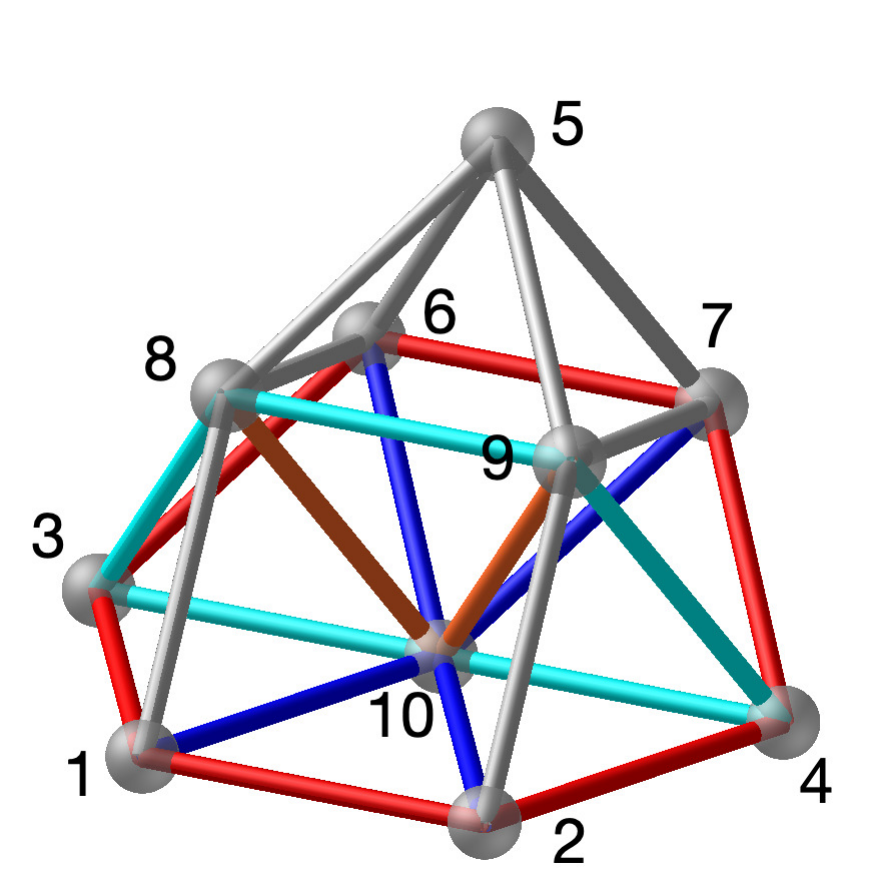}
	}
	\hfil
	\subfloat[]{
		\includegraphics[width=0.2\linewidth]{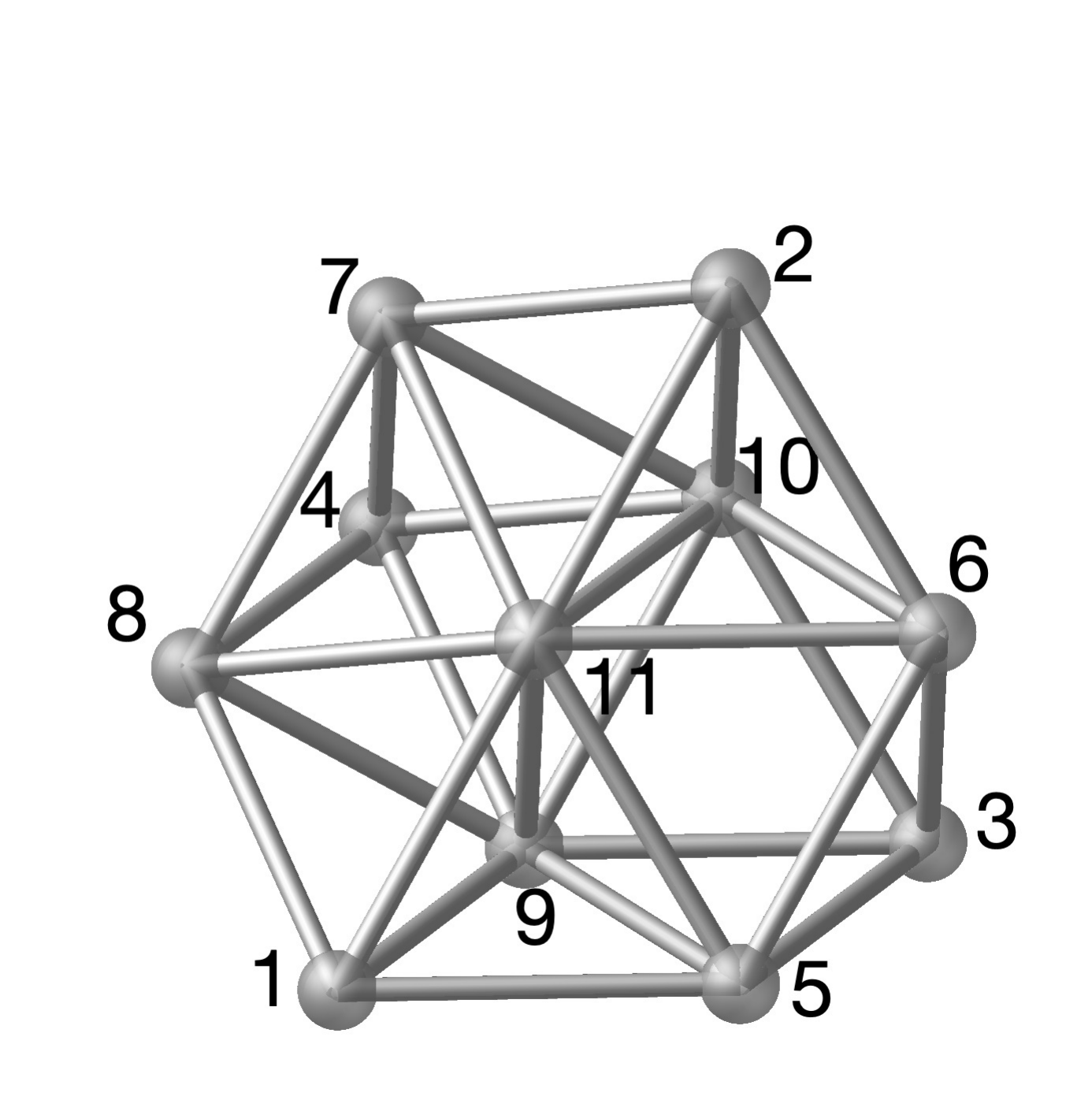}
	}\\
	\subfloat[]{
		\begin{tikzpicture}[scale=0.5]
			
			\coordinate (A) at (0.0, 0.0);
			\coordinate (B) at (1.0, 0.0);
			\coordinate (C) at (-0.15, 0.83);
			\coordinate (D) at (0.79, 1.17);
			\coordinate (E) at (-0.57, 1.56);
			\coordinate (F) at (0.20, 2.20);
			\coordinate (G) at (-1.21, 2.09);
			\coordinate (H) at (-0.72, 2.96);
			\coordinate (I) at (-2.00, 2.38);
			\coordinate (J) at (-1.83, 3.37);
			
			\draw[ultra thick, red!70.71] (A) -- (B); 
			\draw[ultra thick, red!100] (A) -- (C); 
			\draw[ultra thick, red!88.39] (C) -- (D); 
			\draw[ultra thick, red!50.01] (D) -- (B); 
			\draw[ultra thick, blue!70.71] (C) -- (B); 
			\draw[ultra thick, blue!70.71] (A) -- (D); 
			
			\draw[ultra thick, red!25.00] (E) -- (C); 
			\draw[ultra thick, red!88.39] (E) -- (F); 
			\draw[ultra thick, red!12.50] (D) -- (F); 
			\draw[ultra thick, blue!17.68] (C) -- (F); 
			\draw[ultra thick, blue!17.68] (E) -- (D); 
			
			\draw[ultra thick, red!100] (E) -- (G); 
			\draw[ultra thick, red!88.39] (G) -- (H); 
			\draw[ultra thick, red!50] (F) -- (H); 
			\draw[ultra thick, blue!70.71] (E) -- (H); 
			\draw[ultra thick, blue!70.71] (G) -- (F); 
			
			\draw[ultra thick, red!25] (I) -- (G); 
			\draw[ultra thick, blue!17.68] (I) -- (H); 
			\draw[ultra thick, red!12.5] (J) -- (H); 
			\draw[ultra thick, blue!17.68] (G) -- (J); 
			\draw[ultra thick, red!17.68] (J) -- (I); 
			
			\foreach \point in {A, C, E, G, I} {
				\fill[black] (\point) circle (2.5pt);
				\node[left] at (\point) {$\point$};
			}
			\foreach \point in {B, D, F, H, J} {
				\fill[black] (\point) circle (2.5pt);
				\node[right] at (\point) {$\point$};
			}
		\end{tikzpicture}
	}
	\hfil
	\subfloat[]{
		\includegraphics[width=0.2\linewidth]{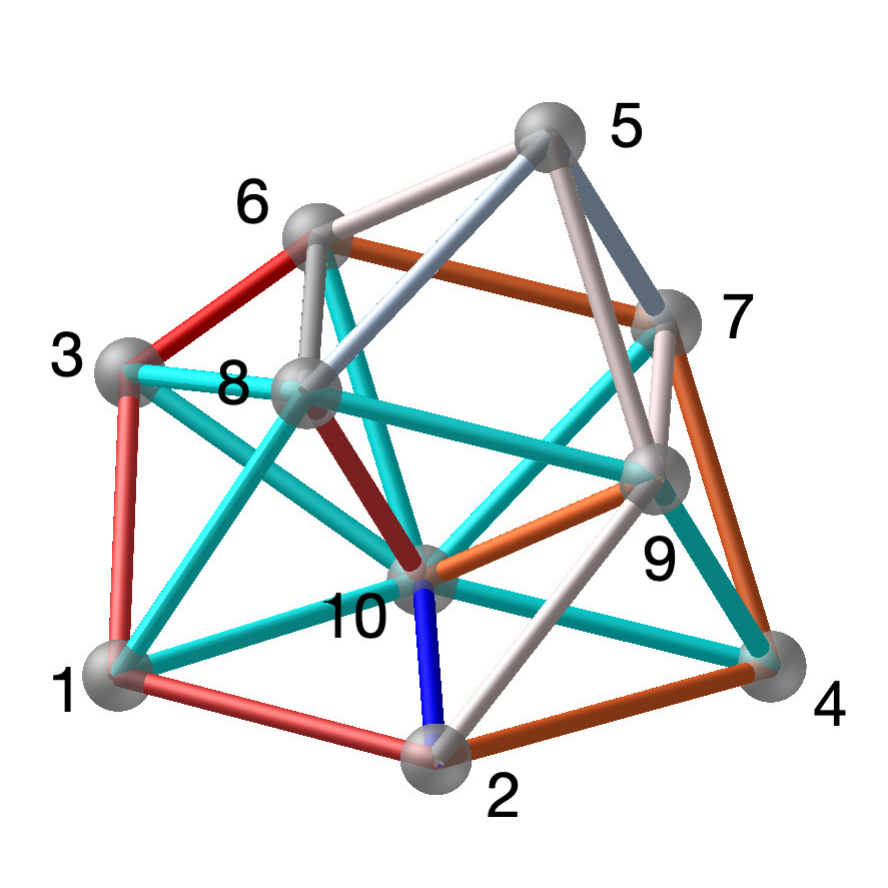}
	}
	\hfil
	\subfloat[]{
		\includegraphics[width=0.2\linewidth]{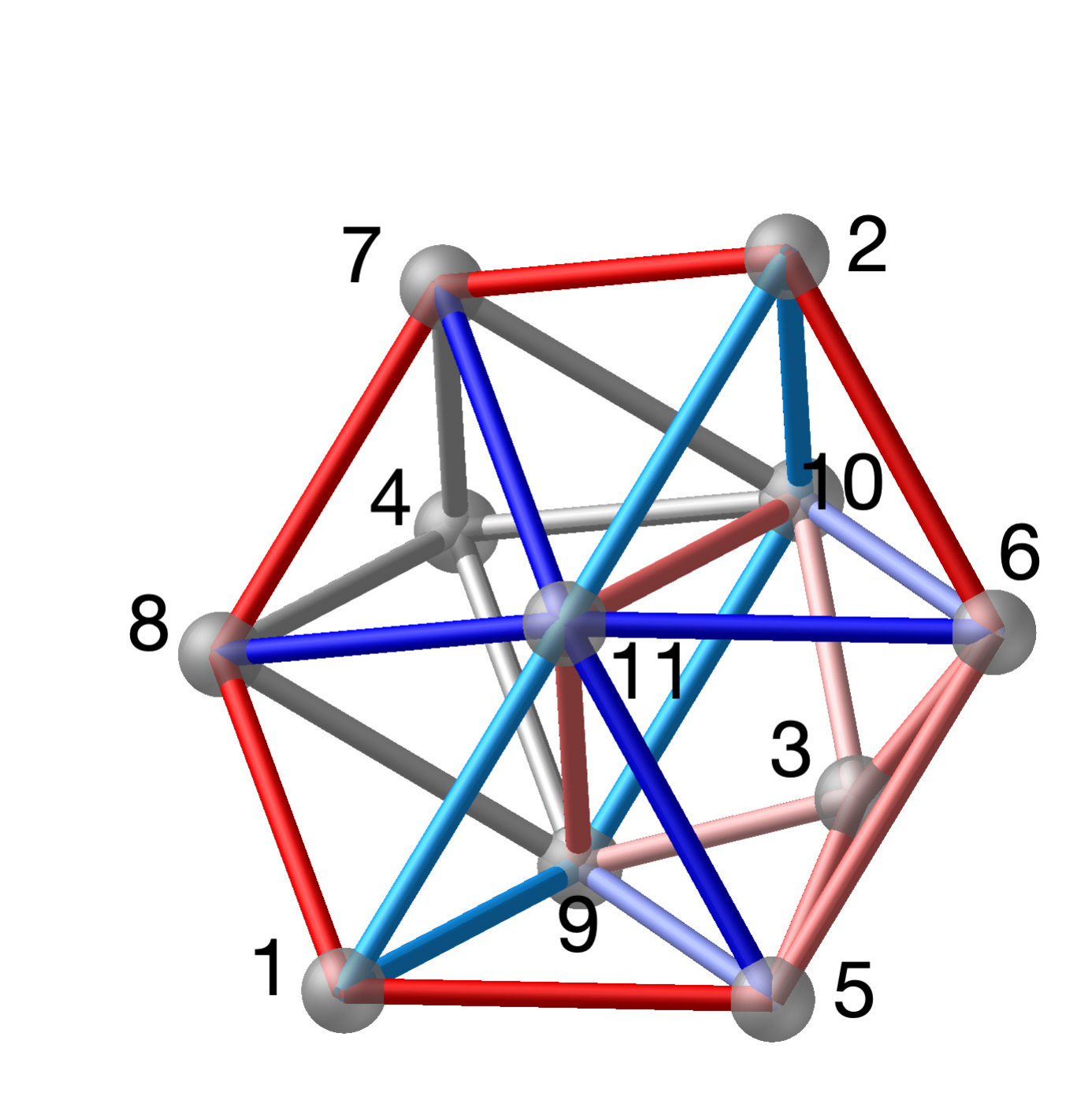}
	}
	\caption{Designing stress ratios on target edges. (a)-(c) are over-constrained  frameworks used as initial conditions for the optimization, and (d)-(f) are solutions to \eqref{eqn:stress-design} with desired self-stresses. The color scheme is the same as the one in \cref{fig:3d-opt-length}.}
	\label{fig:stress-design}
\end{figure}

\paragraph{Examples obtained by stress design} We apply our stress design approach  by solving \eqref{eqn:stress-design} for a 2D stacked square framework and two 3D over-constrained frameworks from the sphere packing dataset\footnote{Figure \ref{fig:stress-design}(b) is the 4th case in $N=10$ and \cref{fig:stress-design}(c) is the 34th case in $N=11$, which can be found in \url{https://personal.math.ubc.ca/~holmescerfon/packings.html}.}, as shown in \cref{fig:stress-design}. The 2D stacked square in \cref{fig:stress-design}(a) has 10 vertices, 21 edges, and a 4-dimensional space of self-stresses. Originally, there is no self-stress with the target ratio 8:4:2:1:8:4:2:1 on edges $AC, BD, CE, DF$, $EG, FH, GI, HJ$ (dashed in \cref{fig:stress-design}(a)). To get a structure with this stress ratio, we set these edges free and solve the optimization problem \eqref{eqn:stress-design}. The resulting optimal framework, shown in \cref{fig:stress-design}(d), has a self-stress matching the specified ratio. (It also has 3 other self-stresses, which are uncontrolled.)

For the 3D framework in \cref{fig:stress-design}(b) with 10 vertices, 25 edges, and a 1-dimensional space of self-stresses, the initial self-stress has a ratio 1:1 on edges 2--4 and 3--6. By applying a target ratio of 1:2 in the optimization, we obtained the framework in \cref{fig:stress-design}(e) with the desired self-stress ratio. 

Lastly, for the 3D framework in \cref{fig:stress-design}(c) with 11 vertices, 29 edges, and a 2-dimensional space of self-stresses, we set a target ratio 2:3:2 for edges 3--5, 2--6, and 3--6. The resulting optimal framework, shown in \cref{fig:stress-design}(f), has a 3-dimensional space of self-stresses, indicating it gained one independent self-stress and a first-order flex. Despite the presence of a first-order flex, the framework remains prestress stable. 

\section{Creating third-order rigid frameworks}\label{sec:higher-rigidity}

The frameworks we have designed so far are all prestress stable, and hence, as discussed in Section \ref{subsec:prelim-rigidity}, they are all second-order rigid. This raises the natural question: how can we design a bar framework that is rigid at an order higher than 2?

We note that the concept of higher-order rigidity has received much attention in the literature, and for a long time there was no consensus about how to properly define it beyond order 2 \cite{connelly1994higher}. We follow the recent theory introduced by \cite{gortler2025higher} which defines orders of rigidity beyond 2 in a way which is robust to previous concerns. 

We present a method to construct isostatic bar frameworks that are so-called third-order rigid, using the definition of this order in \cite{gortler2025higher}.
To do this, we seek a special combination of edge lengths such that $f_1(p)$, the squared length of an edge, reaches a critical point with locally \textbf{cubic growth} up to rigid body motions, while the lengths of all other edges are constrained as in \eqref{eqn:opt}. Under certain conditions this framework can be shown to be third-order rigid. 

\subsection{Theory: a characterization of third-order rigidity via motion on a feasible set}

We first present a condition that guarantees a framework is third-order rigid. In the next section we show how to design frameworks that satisfy this condition.

To simplify the treatment of rigid body motions, we adopt a pinning scheme that fixes $D = d(d+1)/2$ coordinates to zero, via pinning functions $\{g_j\}_{j=1}^D$, as described in Appendix~\ref{app:pinning}. This pinning preserves the order of rigidity \cite{gortler2025higher}. We write the whole collection of length and pinning functions as $\{f_i\}_{i=1}^{m+D}$. We'll consider isostatic graphs, so $m+D=dn$.

First we give a characterization of third-order rigidity. 
%
For this, we must define higher-order flexes. 

\begin{definition}
	A $(1,k)$-flex is a trajectory $p(t)$ such that $p(0) = p$, $p'(0)\neq 0$, and $\frac{d^j}{dt^j}f_i(p(t))\big|_{t=0} = 0$ for $j=1,\ldots,k$, and for all $i=1,\ldots,m+D$.
\end{definition}

From the definition, it is clear that (with pinning constraints), an infinitesimal flex is a $(1,1)$ flex. We write the space of infinitesimal flexes as 
\begin{equation}\label{eqn:pin-1-1-flex}
	K(p) := \{v \in \mathbb{R}^{nd} \:|\: \nabla f_i(p) v = 0, \text{ for all }i=1,\dots,m+D\}.
\end{equation}

\begin{theorem}[\edit{\cite{gortler2025higher}}, Theorems 2.20, 5.3]\label{thm:korder}
	If $\dim K(p) = 1$,  and there exists an integer $k$ such that the pinned framework has a $(1,k-1)$-flex but no $(1,k)$-flex, then the framework has rigidity order $k$. Furthermore, it is rigid. Moreover, if there exists a $(1,k)$ flex, then the framework is either flexible, or it is rigid at a higher order than $k$. 
\end{theorem}

While this result does not define higher-order rigidity, we may take it to be a characterization of higher-order rigidity in the case when $\dim K(p) = 1$. See \cite{gortler2025higher} for a complete characterization of higher-order rigidity.

Our primary result considers  the behaviour of the function $f_1(p(t))$ along a trajectory such that the length and pinning constraints hold exactly:
\begin{equation}\label{eq:con}
	f_2(p(t))-f_2(p(0)) = \cdots = f_{m+D}(p(t))-f_{m+D}(p(0)) = 0. 
\end{equation}
Inspired by the language of constrained optimization, we will say the \emph{LICQ condition holds for the constraint set \eqref{eq:con}} when 
$\{\nabla f_i(p)\}_{i=2}^{m+D}$ are linearly independent. Notice that in this case, the set of solutions to \eqref{eq:con} is locally a one-dimensional manifold, by the implicit function theorem.

\begin{lemma}\label{lem:dimK}
	Consider an isostatic framework $(p,E)$ and suppose the LICQ condition holds for constraint set \eqref{eq:con}. Then $\dim K(p) \leq 1$.
\end{lemma}

\begin{proof}
	This follows from Remark \ref{rmk:fundamental-lin-alg}: if $\dim \nullspace R^T(p) \leq 1$, the fundamental theorem of linear algebra implies that $\dim (\nullspace R(p)\cap \mathcal T(p)^\perp) \leq 1$. 
	But we also have $\dim K(p) \leq \dim (\nullspace R(p)\cap \mathcal T(p)^\perp)$. 
\end{proof}

Our main result is the following. 

\begin{theorem}\label{thm:cubic-growth}
	Consider an isostatic framework $(p^*,E)$ 
	such that the LICQ condition holds for constraint set \eqref{eq:con}.  Suppose there is an analytic path $p(t)$ with $p(0) = p^*$ and $p'(0)\neq 0$, satisfying the constraints \eqref{eq:con}, and such that $f_1(p(t))-f_1(p^*) = a_3 t^3 + o(t^3)$ in a neighborhood of $t=0$. Then, the  framework $(p^*,E)$ is not second-order rigid. Furthermore, if $a_3 \neq 0$, then the framework is third-order rigid.
\end{theorem}

\begin{remark}
	Instead of assuming the LICQ condition holds as a condition of the theorem, one could instead assume that $\mbox{dim}K(p)=1$, plus the additional condition that the unique self-stress $w$ has $w_1\neq 0$.
\end{remark}

\begin{proof}[Proof of Theorem \ref{thm:cubic-growth}]
	First note that by definition, $p(t)$ is a $(1,2)$ flex,  and hence also a $(1,1)$ flex. Therefore, by Lemma \ref{lem:dimK}, $\mbox{dim} \: K(p^*) = 1$, and so by Theorem \ref{thm:korder}, the framework is not second-order rigid.

	Now, we wish to argue that when $a_3\neq 0$, there does not exist a $(1,3)$ flex. To do this, we use Theorem 2.22  of \cite{gortler2025higher} and the discussion preceding it, which imply that, when $\dim K(p^*) = 1$, then if there is a $(1,k)$ flex, then there is a unique $(1,k)$ flex of the form $p(t)= p + p't + p''t^2 + \cdots + p^{(k)}t^k$ with $p^{(j)}\in \bar{K}$ for $j\geq 2$, up to scaling of $p'$. Here $\bar K$ is any space complementary to $K(p^*)$. 
	
	We use this fact by trying to solve for the unique such $(1,3)$ flex, assuming it exists. We start by Taylor-expanding our trajectory as $p(t) = p^* + p' t + p'' t^2 + p''' t^3+\cdots$. 
	We take $\bar K$ to be any space containing $p''$ and which is complementary to $K(p^*)$. 
	Now, we look for a trajectory $q(t) = p^* + p't + p''t^2 + q'''t^3$ which is a $(1,3)$ flex, and with $q'''\in \bar K$. Note that by construction, any such $q(t)$ is already a $(1,2)$ flex.

	Observe that, for the length constraints, 
	\[\frac{d^3}{dt^3}f_i(q(t))\big|_{t=0} = \big\langle p_{i,1}^* - p_{i,2}^*, q_{i,1}''' - q_{i,2}'''\big\rangle + 3\big\langle p'_{i,1} - p'_{i,2}, p_{i,1}'' - p_{i,2}''\big\rangle.
	\]
	For the pinning constraints, we have $\frac{d^3}{dt^3}f_i(q(t))\big|_{t=0} = q_{i'}'''$ if $i'$ is the  coordinate pinned by function $f_i$.
	
	Writing the equations in matrix form, we have that for a $(1,3)$ flex, $q'''$ must solve 
	\begin{equation}\label{eqn:3rd-condition-1}
		\bar R(p^*) q''' + b = 0, 
	\end{equation}
	where $b$ is a vector depending on $p',p''$ (specifically, for the length constraints, $b_i= 3\big\langle p'_{i,1} - p'_{i,2}, p_{i,1}'' - p_{i,2}''\big\rangle$, and for the pinning constraints, $b_i=0$). Here $\bar R(p^*)$ is a matrix whose rows are $\nabla f_i(p^*)$; it is a sort of ``enhanced'' rigidity matrix, with additional rows corresponding to the pinning constraints\footnote{The $\bar R(p^*)$ is in fact $R_g(p^*)$ defined near \eqref{eqn:pinning-condition} in Appendix~\ref{app:pinning}.}.
	
	We are given as conditions of the theorem that 
	$\frac{d^3}{dt^3} f_1(p(t))\big|_{t=0} =a_3 \neq 0$, and
	$\frac{d^3}{dt^3} f_i(p(t)) \big|_{t=0} = 0$, $i=2,\ldots,m$. 
	Thus,
	\begin{equation}\label{eqn:3rd-p-1}
		\bar R(p) p''' + b = a_3 e_1,
	\end{equation}
	where $e_1$ is the first standard basis vector.
	Subtracting \eqref{eqn:3rd-condition-1} and \eqref{eqn:3rd-p-1} gives 
	\begin{equation}\label{eqn:3rd-diff}
		\bar R(p) (p''' - q''') = a_3 e_1.
	\end{equation}

	We wish to show there is no solution to \eqref{eqn:3rd-diff}. First, notice that there exists a self-stress $w$ of the (unpinned) framework such that $w_1\neq 0$. This holds because $p(t)$ is a pinned $(1,1)$-flex; by using Proposition 2.5 in \cite{gortler2025higher}, this implies the existence of an unpinned nontrivial infinitesimal flex. Since the framework is isostatic, Remark \ref{rmk:fundamental-lin-alg} implies there is at least one self-stress $w$. We must have $w_1\neq0$ or else the LICQ condition would be violated.  We may append zeros to $\bar w$, as in $\bar w = (w,0,\ldots 0)$ to get a left null vector of $\bar R(p^*)$. 
	
	Now, take the inner product of $\bar w$ with \eqref{eqn:3rd-diff} to get $0=w_1a_3$. Clearly this is impossible since $w_1\neq 0, a_3\neq 0$. Therefore, no such $q'''$ exists. So, the framework is 3rd order rigid.

\end{proof}

\subsection{Algorithm to construct a third-order rigid framework}\label{sec:third-order-algorithm}

\edit{In general, it is difficult to impose $\frac{d^2}{dt^2}f_1(p(t))=0$ directly, since $t$ parametrizes motion along the one-dimensional constraint manifold and the flex $p(t)$ is unknown in advance.}
\edit{We therefore use an idea from bifurcation theory to find the cubic growth described in Theorem~\ref{thm:cubic-growth} by varying a bar length. For a one-dimensional function depending smoothly on a parameter $\mu$, a local maximum and a local minimum can merge as $\mu$ varies, and at a generic merger the function has cubic growth.}
\edit{Such a merger is known as a \textit{fold catastrophe} (see e.g. Figure 1 and Theorem 4.1 of~\cite{golubitsky1978catastrophe}). At the merger, the first and second derivatives of the squared free edge length along the constraint curve vanish, and the third derivative is nonzero.}

Starting with an isostatic (pinned) framework $(p_0,G)$, we 
assume that $M = \{p \in \mathbb{R}^{nd} \:|\: f_i(p) = f_i(p_0), i=2,\dots,m\}$ is a 1-dimensional smooth manifold parameterized locally by $t \in \mathbb{R}$. We consider a trajectory $p(t)$ on $M$, and monitor the  value of $f_1(p(t))$ as $t$ varies. Typically, $f_1(p(t))$ will have one or more local maximum and local minima; as we showed in Section \ref{sec:our-theorems}, at these optima we obtain a framework that is prestress stable.

To obtain a point with cubic behaviour, we choose another edge length to adjust, say edge 2. Our tuning scheme consists of the following steps:
\begin{enumerate}[(1)]
	\item We first find a local parameter $t \in \mathbb{R}$ that characterizes the constraint manifold $M$. In some special cases, the parameter can be an angle in the plot.
	
	\item We plot $f_1(p)$ as a function of the parameter $t$. The graph $f_1(p(t))$ may have multiple local minima and maxima.
	
	\item We change the length of another edge $E_2$, update this new length in the constraint $f_2(p)=l_2^2$ and plot the graph $f_1(p(t))$ on the new constraint manifold.
	
	\item We stop step (3) until we observe that a local minimum and a local maximum  merge into a critical point, resulting in the cubic growth near the critical point. If the bar framework corresponding to the critical point satisfies the condition $\dim K(p) = 1$, then it is third-order rigid, by Theorem \ref{thm:cubic-growth}. 
\end{enumerate}

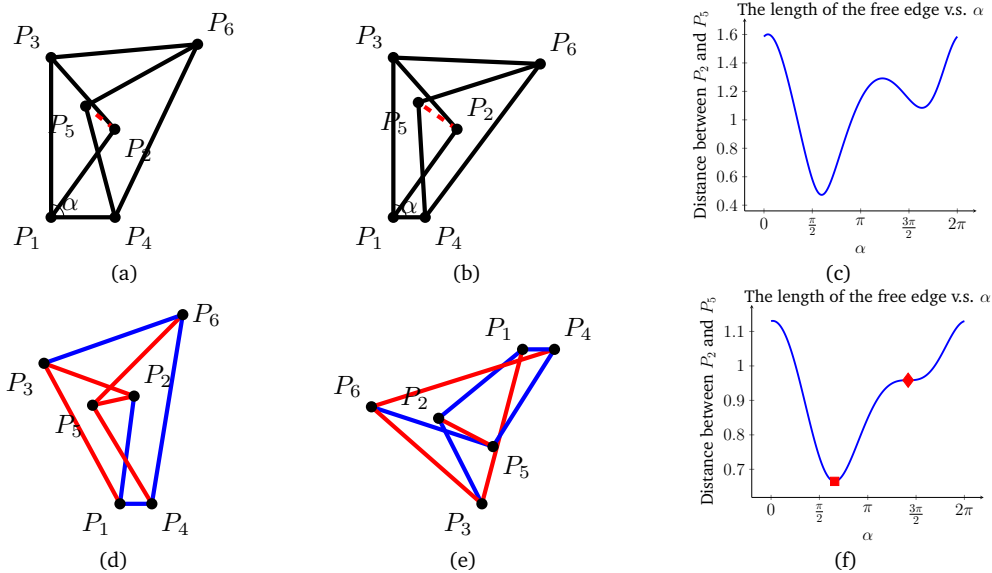
\begin{figure}[!htb]
	\centering
	\captionsetup[subfloat]{farskip=-2pt,captionskip=2pt}
	\subfloat[]{
		\begin{tikzpicture}[scale=0.68]
			\coordinate (A1) at (0,0);
			\coordinate (A2) at (9.9555e-01, 1.3780e+00);
			\coordinate (A3) at (0, 2.5000e+00);
			
			\coordinate (B1) at (1, 0);
			\coordinate (B2) at (5.4036e-01, 1.7403e+00);
			\coordinate (B3) at (2.2906e+00, 2.7082e+00);
			
			\draw[ultra thick,black] (A1) -- (A2); 
			\draw[ultra thick,black] (A2) -- (A3); 
			\draw[ultra thick,black] (A3) -- (A1); 
			
			\draw[ultra thick,black] (B1) -- (B2); 
			\draw[ultra thick,black] (B2) -- (B3);
			\draw[ultra thick,black] (B3) -- (B1);
			
			\draw[ultra thick,red,dashed] (A2) -- (B2); 
			\draw[ultra thick,black] (A1) -- (B1);
			\draw[ultra thick,black] (A3) -- (B3);
			
			\draw (0,0) -- (0.2,0) arc[start angle=0, end angle=90, radius=0.2];
			\node at (0.3, 0.2) {$\alpha$};
			
			\fill[black] (A1) circle (2.5pt);
			\node[below left] at (A1) {$P_1$};
			\fill[black] (A2) circle (2.5pt);
			\node[below right] at (A2) {$P_2$};
			\fill[black] (A3) circle (2.5pt);
			\node[above left] at (A3) {$P_3$};
			
			\fill[black] (B1) circle (2.5pt);
			\node[below right] at (B1) {$P_4$};
			\fill[black] (B2) circle (2.5pt);
			\node[below left] at (B2) {$P_5$};
			\fill[black] (B3) circle (2.5pt);
			\node[above right] at (B3) {$P_6$};
		\end{tikzpicture}
	}
	\hfil
	\subfloat[]{
		\begin{tikzpicture}[scale=0.68]
			\coordinate (A1) at (0,0);
			\coordinate (A2) at (9.9555e-01, 1.3780e+00);
			\coordinate (A3) at (0, 2.5000e+00);
			
			\coordinate (B1) at (4.9815e-01, 0);
			\coordinate (B2) at (3.9104e-01, 1.7968e+00);
			\coordinate (B3) at (2.2978e+00, 2.4002e+00);
			
			\draw[ultra thick,black] (A1) -- (A2); 
			\draw[ultra thick,black] (A2) -- (A3); 
			\draw[ultra thick,black] (A3) -- (A1); 
			
			\draw[ultra thick,black] (B1) -- (B2); 
			\draw[ultra thick,black] (B2) -- (B3);
			\draw[ultra thick,black] (B3) -- (B1);
			
			\draw[ultra thick,red,dashed] (A2) -- (B2); 
			\draw[ultra thick,black] (A1) -- (B1);
			\draw[ultra thick,black] (A3) -- (B3);
			
			\draw (0,0) -- (0.2,0) arc[start angle=0, end angle=90, radius=0.2];
			\node at (0.25, 0.15) {$\alpha$};
			
			\fill[black] (A1) circle (2.5pt);
			\node[below left] at (A1) {$P_1$};
			\fill[black] (A2) circle (2.5pt);
			\node[above right] at (A2) {$P_2$};
			\fill[black] (A3) circle (2.5pt);
			\node[above left] at (A3) {$P_3$};
			
			\fill[black] (B1) circle (2.5pt);
			\node[below right] at (B1) {$P_4$};
			\fill[black] (B2) circle (2.5pt);
			\node[below left] at (B2) {$P_5$};
			\fill[black] (B3) circle (2.5pt);
			\node[above right] at (B3) {$P_6$};
		\end{tikzpicture}
	}
	\hfil
	\subfloat[]{
		\begin{tikzpicture}[scale=0.38]
			\begin{axis}[
				font=\LARGE,
				xlabel={$\alpha$},
				ylabel={Distance between $P_2$ and $P_5$},
				title={The length of the free edge v.s. $\alpha$},
				title style={yshift=-10pt},  
				grid=none,
				axis lines=left,  
				enlargelimits=true,  
				xtick={0, 1.5708, 3.1416, 4.7124, 6.2832},  
				xticklabels={$0$, $\frac{\pi}{2}$, $\pi$, $\frac{3\pi}{2}$, $2\pi$}
				]
				\addplot[mark=none,color=blue, ultra thick] table [col sep=comma, x=X, y=Y] {data-multiple-optima.csv};
			\end{axis}
		\end{tikzpicture}
	}\\
	\subfloat[]{
		\begin{tikzpicture}[scale=0.68]
			\coordinate (A1) at (0,0);
			\coordinate (A2) at (2.2087e-01,1.6856e+00);
			\coordinate (A3) at (-1.1883e+00, 2.1995e+00);
			
			\coordinate (B1) at (4.9815e-01, 0);
			\coordinate (B2) at (-4.2870e-01, 1.5430e+00);
			\coordinate (B3) at (9.8204e-01, 2.9607e+00);
			
			\draw[ultra thick,blue] (A1) -- (A2); 
			\draw[ultra thick,red] (A2) -- (A3); 
			\draw[ultra thick,red] (A3) -- (A1); 
			
			\draw[ultra thick,red] (B1) -- (B2); 
			\draw[ultra thick,red] (B2) -- (B3);
			\draw[ultra thick,blue] (B3) -- (B1);
			
			\draw[ultra thick,red] (A2) -- (B2); 
			\draw[ultra thick,blue] (A1) -- (B1);
			\draw[ultra thick,blue] (A3) -- (B3);
			
			\fill[black] (A1) circle (2.5pt);
			\node[below left] at (A1) {$P_1$};
			\fill[black] (A2) circle (2.5pt);
			\node[above right] at (A2) {$P_2$};
			\fill[black] (A3) circle (2.5pt);
			\node[below left] at (A3) {$P_3$};
			
			\fill[black] (B1) circle (2.5pt);
			\node[below right] at (B1) {$P_4$};
			\fill[black] (B2) circle (2.5pt);
			\node[below left] at (B2) {$P_5$};
			\fill[black] (B3) circle (2.5pt);
			\node[right] at (B3) {$P_6$};
		\end{tikzpicture}
	}
	\hfil
	\subfloat[]{
		\begin{tikzpicture}[scale=0.68]
			\coordinate (A1) at (0,0);
			\coordinate (A2) at (-1.3135e+00, -1.0792e+00);
			\coordinate (A3) at (-6.3624e-01, -2.4177e+00);
			
			\coordinate (B1) at (4.9815e-01, 0);
			\coordinate (B2) at (-4.6325e-01, -1.5217e+00);
			\coordinate (B3) at (-2.3639e+00, -8.9924e-01);
			
			\draw[ultra thick,blue] (A1) -- (A2); 
			\draw[ultra thick,blue] (A2) -- (A3); 
			\draw[ultra thick,red] (A3) -- (A1); 
			
			\draw[ultra thick,blue] (B1) -- (B2); 
			\draw[ultra thick,blue] (B2) -- (B3);
			\draw[ultra thick,red] (B3) -- (B1);
			
			\draw[ultra thick,red] (A2) -- (B2); 
			\draw[ultra thick,blue] (A1) -- (B1);
			\draw[ultra thick,red] (A3) -- (B3);
			
			\fill[black] (A1) circle (2.5pt);
			\node[above left] at (A1) {$P_1$};
			\fill[black] (A2) circle (2.5pt);
			\node[above left] at (A2) {$P_2$};
			\fill[black] (A3) circle (2.5pt);
			\node[below left] at (A3) {$P_3$};
			
			\fill[black] (B1) circle (2.5pt);
			\node[above right] at (B1) {$P_4$};
			\fill[black] (B2) circle (2.5pt);
			\node[below right] at (B2) {$P_5$};
			\fill[black] (B3) circle (2.5pt);
			\node[above left] at (B3) {$P_6$};
		\end{tikzpicture}
	}
	\hfil
	\subfloat[]{
		\begin{tikzpicture}[scale=0.38]
			\begin{axis}[
				font=\LARGE,
				xlabel={$\alpha$},
				ylabel={Distance between $P_2$ and $P_5$},
				title={The length of the free edge v.s. $\alpha$},
				grid=none,
				title style={yshift=-10pt},  
				grid=none,
				axis lines=left,  
				enlargelimits=true,  
				xtick={0, 1.5708, 3.1416, 4.7124, 6.2832},  
				xticklabels={$0$, $\frac{\pi}{2}$, $\pi$, $\frac{3\pi}{2}$, $2\pi$}
				]
				\addplot[mark=none,color=blue, ultra thick] table [col sep=comma, x=X, y=Y] {data-critical.csv};
				
				\addplot[
				only marks,
				mark=diamond*,
				mark size=6pt,
				color=red 
				] coordinates {(4.4550e+00, 9.5854e-01)};
				
				\addplot[
				only marks,
				mark=square*,
				mark size=3.5pt,
				color= red
				] coordinates {(2.06717, 6.6503e-01)};
			\end{axis}
		\end{tikzpicture}
	}
	\caption{(a) An isostatic bar framework. The plot of the free edge, dashed in (a), versus the parameter $\alpha$ is shown in (c); (b) another isostatic bar framework with a different length between $P_1$ and $P_4$. The corresponding plot of the free edge versus $ \alpha$ is shown in (f); (d) a prestress-stable but not first-order bar framework corresponding to the minimum in (f), indicated by a square;  (e) a third-order rigid bar framework corresponding to the critical point in (f), indicated by a diamond.
	}
	\label{fig:higher-rigidity}
\end{figure}

\paragraph{Our third-order rigid framework} We further explain our scheme using a 2D example in \cref{fig:higher-rigidity}. We consider an isostatic bar framework constructed by two triangles $\Delta P_1P_2P_3$ and $\Delta P_4 P_5 P_6$ with three extra edges $P_1 P_4, P_2 P_5$ and $P_3 P_6$. The bar framework is isostatic with 6 vertices and 9 edges. 

We pin $P_1$ at the origin and $P_4$ on the $x$-axis, and choose $P_2 P_5$ as our free edge. The constraint manifold can be parametrized by the angle $\alpha$ marked in \cref{fig:higher-rigidity}. For most of initial bar frameworks, the graph of the free edge as a function of $\alpha$ has a shape plotted in figure of \cref{fig:higher-rigidity}(c) with two local minima and two local maxima.

We now perform step 3. By varying the distance between $P_1$ and $P_4$, the constraint manifold and the graph of the length between $P_2$ and $P_5$ as a function of $\alpha$ will change. We seek a special length between $P_1$ and $P_4$ where a pair of local maximum and minimum annihilate and a critical point with local cubic growth appears (see \cref{fig:higher-rigidity}(f)). We check numerically that at the special framework,  the space of self-stresses and the space of non-trivial first-order flexes are both 1-dimensional. The second-order stress test mentioned in \eqref{eqn:second-order-test} is -7.9219e-05, which is small enough to be considered zero\footnote{The numerical second-order stress test with value around 1e-4 does not vanish because the derivative at the critical point (represented by a diamond in \cref{fig:higher-rigidity})(f) has a numerical error around 1e-8, which cannot be decreased.}. By a continuity argument, there exists a length $E_2$ such that the first and second derivative of $f_1(p(t))$ vanish at a critical point $p^*$ on the constrained manifold. Thus, by applying Theorem \ref{thm:cubic-growth}, there exists a third-order rigid framework close to the one in \cref{fig:higher-rigidity}(e).

\section{Conclusion and open questions}\label{sec:conclusion}

We develop a systematic constrained optimization framework for constructing higher-order rigid structures. Under mild conditions, local solutions to our constrained optimization problem are prestress stable. Our approach bridges rigidity theory and constrained optimization by linking self-stresses with the KKT conditions and prestress stability with the second-order sufficient condition. In addition, we introduce a bifurcation method to construct third-order rigid frameworks.

There are a few open questions to be explored in the future: 
\begin{itemize}
	\item In Section \ref{sec:higher-rigidity}, we present a bifurcation method that constructs third-order rigid frameworks by varying the lengths of two edges. This viewpoint suggests a natural connection to bifurcation theory and catastrophe theory, where varying more edges could give rise to richer bifurcation patterns and potentially yield frameworks with novel rigidity properties.
	
	\item One promising direction is to extend our method to constrained optimization with inequality constraints, thereby enabling the construction of tensegrity structures with higher-order rigidity.
	
	\item \edit{An open problem is to investigate whether our constrained optimization approach can be extended to construct polyhedra whose surfaces are second-order prestress stable as discussed in~\cite{heRobertz2026prestress}.}

	
	\item Most of the higher-order rigid frameworks we constructed have only a one-dimensional space of self-stresses. This arises because the LICQ condition holds at the solutions of our constrained optimization problem (see Lemma \ref{lem:localLICQ}). Our method does not extend to a systematic algorithm for constructing under-constrained or isostatic frameworks with multiple linearly independent self-stresses. Developing such a construction remains an open problem.
\end{itemize}

\edit{Regarding the last point, we briefly explain the difficulty in constructing frameworks with multiple self-stresses, which corresponds to achieving a critical point of multiple constrained optimization problems simultaneously. Specifically, suppose we aim to construct a framework with two linearly independent self-stresses $w, \mu$ that have non-vanishing stresses on the first two edges, i.e. $w_1, w_2, \mu_1, \mu_2 \neq 0$. Without loss of generality, we may set $w_1 = \mu_2 = 1$ and $w_2 = \mu_1 = 0$. According to Theorem~\ref{thm:2nd-sufficient-cond}, such a framework is likely to be a critical point of the following two constrained optimization problems simultaneously:
	\begin{equation}\label{eqn:sim-opt}
		\begin{aligned}
			(P_1) \quad &\min \: f_1(p) & \quad (P_2) \quad &\min \: f_2(p) \\
			&\text{s.t.} \: \: \: f_i(p) = 0, \quad i=3,\dots,m & &\text{s.t.} \: \:  \: f_i(p) = 0, \quad i=3,\dots,m.
		\end{aligned}
	\end{equation}
	However, finding a configuration that is a critical point for multiple optimization problems is highly nontrivial -- both analytically and numerically -- and may require additional symmetry or special structural conditions.}

\begin{figure}[!htb]
	\centering
	\subfloat[]{
		\begin{tikzpicture}[scale=0.16]
			\coordinate (A) at (0, 0);       
			\coordinate (B) at (8, 0);       
			\coordinate (C) at (5, 8);
			\coordinate (D) at (2.8349, 6.75);
			\coordinate (E) at (4.5422, 4.4942);
			\coordinate (F) at (5.4019, 1.5);       
			
			\draw[blue] (5,8) -- (6.5,8) arc[start angle=0, end angle=-95, radius=1.5];
			\node at (7,6.5) {$\theta_1$};
			
			\draw[blue] (8,0) -- (9.5,0) arc[start angle=0, end angle=145, radius=1.5];
			\node at (10,1.5) {$\theta_2$};
			
			\draw[ultra thick, black] (A) -- (B);
			\draw[ultra thick, black] (B) -- (C);
			\draw[ultra thick, black] (C) -- (A);
			
			\draw[ultra thick, red, dashed] (D) -- (A);
			
			\draw[ultra thick, black] (C) -- (D);
			\draw[ultra thick, black] (D) -- (E);
			\draw[ultra thick, black] (C) -- (E);
			
			\draw[ultra thick, red, dashed] (E) -- (F);
			\draw[ultra thick, black] (B) -- (F);
			
			\coordinate (CircleCenter1) at (5, 8);
			
			\draw[blue, thick, dotted] (CircleCenter1) circle (3.5355);
			
			\coordinate (CircleCenter2) at (8, 0);
			
			\draw[blue, thick, dotted] (CircleCenter2) circle (3);
			
			\filldraw[red] (A) circle (9pt) node[black, below] {A};
			\filldraw[red] (B) circle (9pt) node[black, below] {B};
			\filldraw[red] (C) circle (9pt) node[black, above] {C};
			\filldraw[blue] (D) circle (9pt) node[black, left] {D};
			\filldraw[blue] (E) circle (9pt) node[black, below left] {E};
			\filldraw[blue] (F) circle (9pt) node[black, left] {F};
		\end{tikzpicture}
	}\hfil
	\subfloat[]{
		\begin{tikzpicture}[scale=0.16]
			
			\coordinate (A) at (0, 0);       
			\coordinate (B) at (8, 0);       
			\coordinate (C) at (5, 8);
			\coordinate (D) at (3.68,5.88);
			\coordinate (E) at (6.24,4.69);
			\coordinate (F) at (6.95,2.81);       
			
			\draw[ultra thick, black] (A) -- (B);
			\draw[ultra thick, black] (B) -- (C);
			\draw[ultra thick, black] (C) -- (A);
			
			\draw[ultra thick] (D) -- (A);
			
			\draw[ultra thick, black] (C) -- (D);
			\draw[ultra thick, black] (D) -- (E);
			\draw[ultra thick, black] (C) -- (E);
			
			\draw[ultra thick, black] (E) -- (F);
			\draw[ultra thick, black] (B) -- (F);
			
			\coordinate (CircleCenter1) at (5, 8);
			
			\draw[blue, thick, dotted] (CircleCenter1) circle (3.5355);
			
			\coordinate (CircleCenter2) at (8, 0);
			
			\draw[blue, thick, dotted] (CircleCenter2) circle (3);
			
			\filldraw[red] (A) circle (9pt) node[black, below] {A};
			\filldraw[red] (B) circle (9pt) node[black, below] {B};
			\filldraw[red] (C) circle (9pt) node[black, above] {C};
			\filldraw[blue] (D) circle (9pt) node[black, left] {D};
			\filldraw[blue] (E) circle (9pt) node[black, right] {E};
			\filldraw[blue] (F) circle (9pt) node[black, left] {F};
		\end{tikzpicture}
	}\hfil
	\subfloat[]{
		\begin{tikzpicture}[scale=0.16]
			
			\coordinate (A) at (0, 0);       
			\coordinate (B) at (8, 0);       
			\coordinate (C) at (5, 8);
			\coordinate (D) at (6.32,10.12);
			\coordinate (E) at (3.76,11.31);
			\coordinate (F) at (9.05,-2.81);       
			
			\draw[ultra thick, black] (A) -- (B);
			\draw[ultra thick, black] (B) -- (C);
			\draw[ultra thick, black] (C) -- (A);
			
			\draw[ultra thick, black] (C) -- (D);
			\draw[ultra thick, black] (D) -- (E);
			\draw[ultra thick, black] (C) -- (E);
			
			\draw[ultra thick, black] (B) -- (F);
			\draw[ultra thick] (E) -- (F);
			\draw[ultra thick] (D) -- (A);
			
			\coordinate (CircleCenter1) at (5, 8);
			
			\draw[blue, thick, dotted] (CircleCenter1) circle (3.5355);
			
			\coordinate (CircleCenter2) at (8, 0);
			
			\draw[blue, thick, dotted] (CircleCenter2) circle (3);
			
			\filldraw[red] (A) circle (9pt) node[black, left] {A};
			\filldraw[red] (B) circle (9pt) node[black, right] {B};
			\filldraw[red] (C) circle (9pt) node[black, right] {C};
			\filldraw[blue] (D) circle (9pt) node[black, right] {D};
			\filldraw[blue] (E) circle (9pt) node[black, left] {E};
			\filldraw[blue] (F) circle (9pt) node[black, right] {F};
		\end{tikzpicture}
	}\\
	\subfloat[]{
		\begin{tikzpicture}[scale=0.16]
			
			\coordinate (A) at (0, 0);       
			\coordinate (B) at (8, 0);       
			\coordinate (C) at (5, 8);
			\coordinate (D) at (3.68,5.88);
			\coordinate (E) at (6.24,4.69);
			\coordinate (F) at (9.05,-2.81);       
			
			\draw[ultra thick, black] (A) -- (B);
			\draw[ultra thick, black] (B) -- (C);
			\draw[ultra thick, black] (C) -- (A);
			
			\draw[ultra thick] (D) -- (A);
			
			\draw[ultra thick, black] (C) -- (D);
			\draw[ultra thick, black] (D) -- (E);
			\draw[ultra thick, black] (C) -- (E);
			
			\draw[ultra thick, black] (E) -- (F);
			\draw[ultra thick, black] (B) -- (F);
			
			\coordinate (CircleCenter1) at (5, 8);
			
			\draw[blue, thick, dotted] (CircleCenter1) circle (3.5355);
			
			\coordinate (CircleCenter2) at (8, 0);
			
			\draw[blue, thick, dotted] (CircleCenter2) circle (3);
			
			\filldraw[red] (A) circle (9pt) node[black, left] {A};
			\filldraw[red] (B) circle (9pt) node[black, right] {B};
			\filldraw[red] (C) circle (9pt) node[black, above] {C};
			\filldraw[blue] (D) circle (9pt) node[black, left] {D};
			\filldraw[blue] (E) circle (9pt) node[black, right] {E};
			\filldraw[blue] (F) circle (9pt) node[black, right] {F};
		\end{tikzpicture}
	}
	\hfil
	\subfloat[]{
		\begin{tikzpicture}[scale=0.16]
			
			\coordinate (A) at (0, 0);       
			\coordinate (B) at (8, 0);       
			\coordinate (C) at (5, 8);
			\coordinate (D) at (6.32,10.12);
			\coordinate (E) at (3.76,11.31);
			\coordinate (F) at (6.95,2.81);       
			
			\draw[ultra thick, black] (A) -- (B);
			\draw[ultra thick, black] (B) -- (C);
			\draw[ultra thick, black] (C) -- (A);
			
			\draw[ultra thick, black] (C) -- (D);
			\draw[ultra thick, black] (D) -- (E);
			\draw[ultra thick, black] (C) -- (E);
			
			\draw[ultra thick, black] (B) -- (F);
			\draw[ultra thick] (E) -- (F);
			\draw[ultra thick] (D) -- (A);
			
			\coordinate (CircleCenter1) at (5, 8);
			
			\draw[blue, thick, dotted] (CircleCenter1) circle (3.5355);
			
			\coordinate (CircleCenter2) at (8, 0);
			
			\draw[blue, thick, dotted] (CircleCenter2) circle (3);
			
			\filldraw[red] (A) circle (9pt) node[black, left] {A};
			\filldraw[red] (B) circle (9pt) node[black, right] {B};
			\filldraw[red] (C) circle (9pt) node[black, right] {C};
			\filldraw[blue] (D) circle (9pt) node[black, right] {D};
			\filldraw[blue] (E) circle (9pt) node[black, left] {E};
			\filldraw[blue] (F) circle (9pt) node[black, right] {F};
		\end{tikzpicture}
	}
	\hfil
	\subfloat[]{
		\begin{tikzpicture}[scale=0.16]
			\coordinate (A) at (0, 0);       
			\coordinate (B) at (8, 0);       
			\coordinate (C) at (5, 8);
			\coordinate (D) at (2.8349, 6.75);
			\coordinate (E) at (10-4.5422, 4.4942);
			\coordinate (F) at (5.4019, 1.5);       
			
			\draw[blue] (5,8) -- (6.5,8) arc[start angle=0, end angle=-85, radius=1.5];
			\node at (7,6.5) {$\theta_1$};
			
			\draw[blue] (8,0) -- (9.5,0) arc[start angle=0, end angle=145, radius=1.5];
			\node at (10,1.5) {$\theta_2$};
			
			\draw[ultra thick, black] (A) -- (B);
			\draw[ultra thick, black] (B) -- (C);
			\draw[ultra thick, black] (C) -- (A);
			
			\draw[ultra thick, red, dashed] (D) -- (A);
			
			\draw[ultra thick, black] (C) -- (D);
			\draw[ultra thick, black] (D) -- (E);
			\draw[ultra thick, black] (C) -- (E);
			
			\draw[ultra thick, red, dashed] (E) -- (F);
			\draw[ultra thick, black] (B) -- (F);
			
			\coordinate (CircleCenter1) at (5, 8);
			
			\draw[blue, thick, dotted] (CircleCenter1) circle (3.5355);
			
			\coordinate (CircleCenter2) at (8, 0);
			
			\draw[blue, thick, dotted] (CircleCenter2) circle (3);
			
			\filldraw[red] (A) circle (9pt) node[black, below] {A};
			\filldraw[red] (B) circle (9pt) node[black, below] {B};
			\filldraw[red] (C) circle (9pt) node[black, above] {C};
			\filldraw[blue] (D) circle (9pt) node[black, left] {D};
			\filldraw[blue] (E) circle (9pt) node[black, below left] {E};
			\filldraw[blue] (F) circle (9pt) node[black, left] {F};
		\end{tikzpicture}
	}
	\caption{\edit{Frameworks with two self-stresses: (a) the initial framework with $\angle DCE = \angle ACB$ where the critical points of the two constrained optimizations can be achieved simultaneously; (b--e): frameworks with two self-stresses. Examples (b,c) are rigid while examples (d,e) are flexible; (f) another initial framework with $\angle DCE \neq \angle ACB$ where the critical points of the two constrained optimizations cannot be achieved simultaneously.}
	}
	\label{fig:tangent-form}
\end{figure}
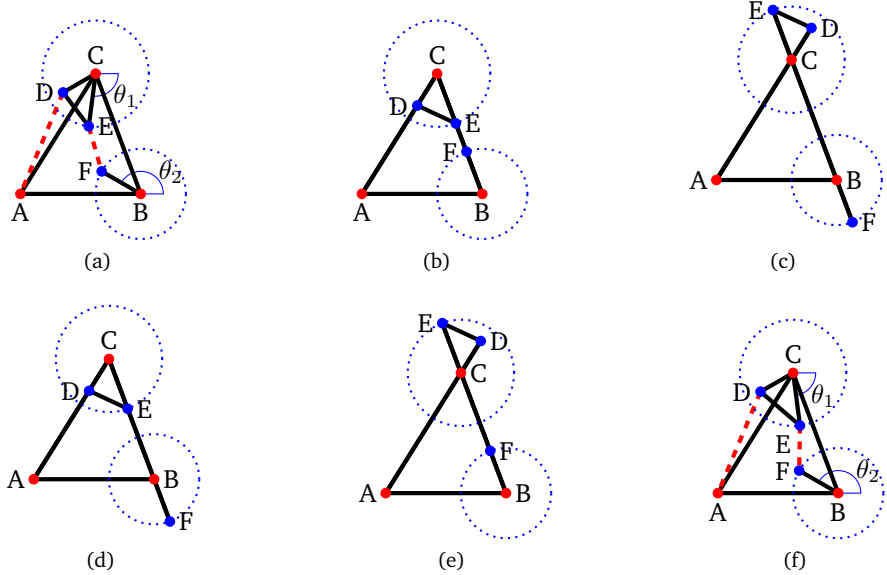

\edit{We illustrate this difficulty with a detailed 2D example in \cref{fig:tangent-form}. Starting from the framework shown in \cref{fig:tangent-form}(a), an isostatic structure with six vertices and nine edges, we aim to construct a framework with two self-stresses using \eqref{eqn:sim-opt}. We deliberately choose the angle $\angle ACB$ equal to $\angle DCE$ with a reason to be explained shortly. We allow two edges $AD,EF$ to change their lengths and choose  $f_1(p),f_2(p)$ to represent the squared lengths of the edges $AD$ and $EF$. Since $\angle ACB = \angle DCE$, $AD$ and $EF$ achieve the smallest lengths at the same structure shown in \cref{fig:tangent-form}(b), which in fact correspond to a rigid framework with two linearly independent self-stresses. However, if we start with a structure $\angle ACB \neq \angle DCE$, the two constrained optimizations in \eqref{eqn:sim-opt} cannot achieve their minimum simultaneously, which partly explains the requirement of special geometry when we use \eqref{eqn:sim-opt} to construct frameworks with multiple self-stresses.
	
	It turns out that we can replace the minimum in \eqref{eqn:sim-opt} to be any kind of critical point -- minimum, maximum, or saddle point. When we search for critical points for $f_1(p)$ and $f_2(p)$ on the constrained manifold with $f_i(p) = 0$ for $i=3,\dots,m$, there are three more solutions as critical points of $f_1(p)$ and $f_2(p)$ at the same time. Hence, there are three more frameworks two linearly independent self-stresses as shown in \cref{fig:tangent-form}(c-e). In fact, the framework in \cref{fig:tangent-form}(c) is rigid and achieves maximum for both $P_1$ and $P_2$ in \eqref{eqn:sim-opt}, while the frameworks in \cref{fig:tangent-form}(d-e) are flexible\footnote{The frameworks in \cref{fig:tangent-form}(d-e) resemble four-bar linkages and their flexes correspond to the two branches of motion in a four-bar linkage.} and achieve a minimum or a maximum for $P_1$ and a saddle point for $P_2$ in \eqref{eqn:sim-opt}. 
}

\appendix
\renewcommand{\theequation}{\thesection.\arabic{equation}}

\section{Proof of Theorem \ref{thm:secondordersufficient}}\label{app:proof-thm216}

Our goal is to show that any pair of $(p^*, \lambda^*)$ satisfying \eqref{eqn:secondsufficient} and the KKT condition \eqref{eqn:KKT} must be a strict local solution up to rigid body motions for $P_k^+$ in \eqref{eqn:opt} (the statement for $P_k^-$ can be proved similarly). We review the condition \eqref{eqn:secondsufficient} by providing an equivalent condition as follows
\begin{equation*}
	v^T \nabla^2_p L(p^*,\lambda^*) v \geq \sigma \|v\|^2,\qquad \forall \:v \in C_k(p^*),
\end{equation*}
with some constant $\sigma > 0$. 

Now we start proving Theorem \ref{thm:secondordersufficient} by contradiction: we assume that there exists a sequence $p_n \rightarrow p^*$ such that (1) $f_k(p_n) \leq f_k(p^*)$; (2) $f_i(p_n) = f_i(p^*)$ with $i\neq k$; and (3) $p_n$ are not congruent to $p^*$. We modify $p_n$ to obtain a new sequence $\widehat{p}_n$ as follows
\begin{equation*}
	\widehat{p}_n = \underset{\widetilde{p}\text{ congruent to }p_n}{\text{argmin}}\|\widetilde{p}-p^*\|.
\end{equation*}
Consequently, the new sequence $\widehat{p}_n
$ satisfies (a) $\widehat{p}_n - p^* \perp \mathcal{T}(\widehat{p}_n)$, where $\mathcal{T}(\widehat{p}_n)$ is the subspace of trivial flexes at $\widehat{p}_n$; (b) $f_i(\widehat{p}_n) = f_i(p_n) = f_i(p^*)$ for $i\neq k$ and $f_k(p_n) = f_k(\widehat{p}_n)$; (c) $\widehat{p}_n$ also converges to $p^*$, since $\|\widehat{p}_n - p^*\| \leq \|p_n - p^*\|$. 

We observe that the limiting directional vector $d = \lim_{n \rightarrow \infty} \frac{\widehat{p}_n - p^*}{\|\widehat{p}_n - p^*\|}$ (by passing to a subsequence) satisfies $d \perp \nabla f_i(p^*)$ for $i\neq k$. To see why, we use the identity $f_i(\widehat{p}_n) = f_i(p^*)$ in (b) and apply a Taylor series expansion, which gives
\begin{equation*}
	0 = \nabla f_i(p^*) (\widehat{p}_n-p^*) + \frac{1}{2} (\widehat{p}_n-p^*)^T \nabla^2 f_i(p^*) (\widehat{p}_n-p^*) + o(\|\widehat{p}_n - p^*\|^2).
\end{equation*}
We obtain $d \perp \nabla f_i(p^*)$ by dividing both sides by $\|\widehat{p}_n - p^*\|$ and taking $n \rightarrow \infty$. 

Now we form a contradiction by showing $f_k(p_n) > f_k(p^*)$ for sufficiently large $n$. We use the fact that $\lambda^*$ is a self-stress and $f_i(p_n) = f_i(\widehat{p}_n)=f_i(p^*)$ for $i\neq k$ and obtain
\begin{equation*}
	\begin{aligned}
		f_k(p_n) - f_k(p^*) &= f_k(\widehat{p}_n)-f_k(p^*) = L(\widehat{p}_n,\lambda^*) - L(p^*,\lambda^*) \\
		&= \frac{1}{2} (\widehat{p}_n-p^*)^T \nabla^2 L(p^*) (\widehat{p}_n-p^*) + o(\|\widehat{p}_n - p^*\|^2).
	\end{aligned}
\end{equation*}
It is sufficient to show that
\begin{equation}\label{eqn:cst-D}
	D:=d^T \nabla^2 L(p^*, \lambda^*) d = \lim_{k \rightarrow \infty} \frac{(\widehat{p}_k-p^*)^T \nabla^2 L(p^*) (\widehat{p}_k-p^*)}{\|\widehat{p}_k - p^*\|^2} > 0.
\end{equation}
Using \eqref{eqn:cst-D}, we conclude our proof of Theorem \ref{thm:secondordersufficient} since 
\begin{equation*}
	f_k(\widehat{p}_n) - f_k(p^*) > D/4\|\widehat{p}_n - p^*\|^2
\end{equation*}
holds for sufficiently large $n$, which violates our assumption in (1).

To prove \eqref{eqn:cst-D}, the key is to show that $d \notin \mathcal{T}(p^*)$ otherwise $D$ automatically vanishes. Since $d \perp \nabla f_i(p^*)$ with $i\neq k$, we have $d \in C_k(p^*) \oplus \mathcal{T}(p^*)$. Therefore, we can separate $d$ as $d = v + r$ with $v \in C_k(p^*)$ and $r \in \mathcal{T}(p^*)$. A standard calculation yields (the proof is postponed at the end of the proof)
\begin{equation}\label{eqn:lagrangian-rigid-body}
	v^T \nabla^2 L(p^*,\lambda^*) v = (v+r)^T \nabla^2 L(p^*,\lambda^*) (v+r).
\end{equation}
Therefore, we only need to show that $v \neq 0$, which is equivalent to show $d \notin \mathcal{T}(p^*)$. 

To prove $d \notin \mathcal{T}(p^*)$, we need the projection operator $P_{\mathcal{T}(p)}$ which maps any vector to its projected component in $\mathcal{T}(p)$. Since $P_{\mathcal{T}(p)}:\mathbb{R}^{nd} \rightarrow \mathcal{T}(p)$ is continuous with respect to $p$, we have 
\begin{equation}\label{eqn:proj-bd}
	\|P_{\mathcal{T}(\widehat{p}_n)} - P_{\mathcal{T}(p^*)}\| < 1/2
\end{equation}
for sufficiently large $n$. 

Now we prove $d \notin \mathcal{T}(p^*)$ by contradiction and assume $d \in \mathcal{T}(p^*)$. We first observe that
\begin{equation*}
	1=\|P_{\mathcal{T}(p^*)} d\| =\|P_{\mathcal{T}(p^*)} d-P_{\mathcal{T}(\widehat{p}_n)} d_n\|, 
\end{equation*}
where $d_n = \frac{\widehat{p}_n - p^*}{\|\widehat{p}_n - p^*\|}$. The last equality holds since $\widehat{p}_n - p^* \perp \mathcal{T}(\widehat{p}_n)$. Therefore, we have
\begin{equation*}
	\begin{aligned}
		1=\|P_{\mathcal{T}(p^*)} d-P_{\mathcal{T}(\widehat{p}_n)} d_n\| &\leq \|P_{\mathcal{T}(p^*)} d-P_{\mathcal{T}(p^*)} d_n\| + \|P_{\mathcal{T}(p^*)} d_n-P_{\mathcal{T}(\widehat{p}_n)} d_n\|\\
		&\leq \|d-d_n\| + \|P_{\mathcal{T}(\widehat{p}_n)} - P_{\mathcal{T}(p^*)}\|,
	\end{aligned}
\end{equation*}
where the limit of the upper bound is less than $1/2$ due to \eqref{eqn:proj-bd}. Thus, we have $d \notin \mathcal{T}(p^*)$.


Lastly, we prove \eqref{eqn:lagrangian-rigid-body} by taking a smooth deformation $q(t)$ such that $q(0) = p^*, q '(0)=v$ and another path $q_r(t)$ with $q_r(0) = p^*, q_r'(0) = v+r$ such that $q(t),q_r(t)$ are congruent, i.e. 
\begin{equation}\label{eqn:congruent-motion}
	f_i(q(t)) = f_i(q_r(t)), \qquad i=1,\dots,m.
\end{equation}
Let us briefly explain why a congruent path $q_r(t)$ can be constructed from $q(t)$. Since $r \in \mathcal{T}(p^*)$, we decompose $r = r_t + r_r$, where $r_t$ is the translational component and $r_r = \Omega p^*$ is the rotational component. Here $\Omega \in \mathbb{R}^{nd \times nd}$ is block diagonal with skew-symmetric blocks $W \in \mathbb{R}^{d \times d}$ whose associated rotation is given by $R(t) = \exp(Wt) \in \mathbb{R}^{d \times d}$. We observe that $R(0)=I$ and $R'(0)=W$. Thus the congruent path $q_r(t)$ takes the form
\begin{equation*}
	q_r(t) = \Omega(t)\, q(t) + t\,r_t,
\end{equation*}
where $\Omega(t) \in \mathbb{R}^{nd \times nd}$ is block diagonal with diagonal blocks $R(t)$.

By taking $\frac{d^2}{dt^2}|_{t=0}$ on both sides of \eqref{eqn:congruent-motion}, we obtain
\begin{equation}\label{eqn:2nd-inf-rigid}
	\nabla f_i(p^*) q''(0) + v^T \nabla^2 f_i(p^*) v = \nabla f_i(p^*) q_r''(0) + (v+r)^T \nabla^2 f_i(p^*) (v+r).
\end{equation}
Now we multiply each equation by $\lambda_i$ and sum over $i$ to form $\nabla^2 L$. Since $\lambda$ is a self-stress, $\sum_{i} \lambda_i \nabla f_i = 0$, the acceleration term $q_r''(0)$ drops and we achieve \eqref{eqn:2nd-inf-rigid}.

\section{The pinning scheme} \label{app:pinning}

For a framework in $\mathbb{R}^d$ whose vertices have a $d$-dimensional affine span, we pin $D = d(d+1)/2$ coordinates at zero, by setting $p_1$ to the origin, and then for each $2\leq i\leq d$, setting $p_i$ to be in the span of the first $i-1$ coordinate vectors.  For example, for $d=2$, we fix one vertex at the origin and set the $y$-coordinate of another vertex to zero. For $d=3$, we fix six coordinates: place one vertex at the origin, set the $y,z$-coordinates of a second vertex to zero, and the $z$-coordinate of a third vertex to zero.

The constrained optimization problem with pinning is:
\begin{equation}\label{eqn:opt-pinning}
	\begin{aligned}
		(\Pkpm) \qquad  \min_{p \in \mathbb{R}^{nd}}  \pm f_k(p) \qquad
		\text{s.t } \qquad  f_i(p)&=l_i^2, \quad i = 1,\dots,m, \quad i\neq k,\\
		g_j(p) &= 0, \quad j=1,\dots,D,
	\end{aligned}
\end{equation}
where $g_j(p)=0$ are the linear constraints that pin a coordinate at zero. In general, the pinning does not affect the order of rigidity provided $p_1,\ldots,p_{d+1}$ are affinely independent \cite{gortler2025higher}.\footnote{If the vertices have an $l$-dimensional affine span, then fix a smaller number of vertices. See \cite{gortler2025higher} for details.}

We aim to link the second-order sufficient conditions for \eqref{eqn:opt-pinning}, to the conditions for prestress stability. 
We assume the following condition holds:
\begin{equation}\label{eqn:pinning-condition}
	\text{rank }G(p) (T(p))^T = D,
\end{equation}
where
\[
G(p) = \begin{pmatrix}
	\nabla g_1(p) \\\vdots \\ \nabla g_{D}(p)
\end{pmatrix},\qquad
T(p) = \begin{pmatrix}
	t_1(p) \\\vdots \\ t_D(p)
\end{pmatrix}.
\]
	Notice that condition \eqref{eqn:pinning-condition} can only hold when $ \text{rank} \: T(p) = D$, hence when $\dim \mathcal{T}(p) =D$. Thus, we do not consider frameworks where the trivial motions are degenerate, such as vertices on a line in $\R^3$.

	\begin{remark}[The pinning condition viewed as a complementary condition]
		When $\mathcal T(p)$ has full rank, condition \eqref{eqn:pinning-condition} is equivalent  to requiring that $\nullspace \: G(p)$ be complementary to $\mathcal{T}(p)$, as we briefly explain. For the matrix $G(p)$ satisfying \eqref{eqn:pinning-condition}, we observe that $\dim \: G(p) = \dim \mathcal{T}(p) = D$ and $\dim \nullspace \: G(p) = nd-D$. If there exists $v \in \nullspace \: G(p) \cap \mathcal{T}(p)$, then $v \in \nullspace \: G(p) \cap \text{Range} \: (T(p))^T$. 
		Using \eqref{eqn:pinning-condition}, we have
		\begin{equation*}
			v = (T(p))^T a \quad \Rightarrow\quad  0 = G(p) v = G(p)(T(p))^T a \quad \Rightarrow \quad a = 0 \quad \Rightarrow \quad v = 0,
		\end{equation*}
		which indicates the equivalence of \eqref{eqn:pinning-condition} and $\nullspace \: G(p)$ being complementary to $\mathcal{T}(p)$.
	\end{remark}

	For a local solution of $p^*$ of $\Pkpm$ in \eqref{eqn:opt-pinning} satisfying the pinning condition \eqref{eqn:pinning-condition}, the KKT condition and second-order sufficient conditions of \eqref{eqn:opt-pinning} are equivalent to those of the unpinned problem \eqref{eqn:opt}, which are stated as:

	\begin{lemma}\label{lemma:pinning-kkt}
		For a local solution of $p^*$ of $\Pkpm$ in \eqref{eqn:opt-pinning} that satisfies the pinning condition \eqref{eqn:pinning-condition}, the Lagrange multiplier vanishes on the pinning constraints, i.e. any left null vector of $R_g(p)$ in \eqref{eqn:pinning-matrics} with $w \in \mathbb{R}^m, s \in \mathbb{R}^{D}$ and $w^T R(p) + s^T G(p) = 0$ must have $s=0$ and $w \in \nullspace \: (R(p))^T$.
	\end{lemma}

	\begin{lemma}\label{lemma:2nd-sufficient-pin}
		Suppose $p^*$ is a local solution for one of $P_k^\pm$ in \eqref{eqn:opt-pinning}. If the pinning condition \eqref{eqn:pinning-condition} and the KKT condition \eqref{eqn:KKT} hold at $p^*$, the second-order sufficient condition \eqref{eqn:secondsufficient} of the unpinned problem holds if and only if the second-order sufficient condition of \eqref{eqn:opt-pinning} holds, i.e.
		\begin{equation}\label{eqn:2nd-test-pinning}
			v^T \Big(\sum_{i=1}^m \lambda^*_i \nabla^2 f_i(p^*)\Big) v > 0, \qquad \forall \: v \in C_k^\text{pin}(p^*) \text{ and }v\neq 0,
		\end{equation}
		where the pinned critical cone is
		\begin{equation}\label{eqn:critical-cone-pinning}
			C_k^\text{pin}(p^*) = \big\{v \: |\: \nabla f_i(p^*) v = 0, \: \forall i\neq k \text{ and }\nabla g_j(p^*) v = 0, \: \forall j=1,\dots,D\big\}.
		\end{equation}
	\end{lemma}

	The proof Lemma \ref{lemma:pinning-kkt} and \ref{lemma:2nd-sufficient-pin} will be provided shortly. Our main result of this section about the prestress stability of local solutions to \eqref{eqn:opt-pinning} is stated as:

	\begin{proposition}\label{prop:prestress-pin}
		For a local solution $p^*$ for one of $P_k^\pm$ in \eqref{eqn:opt-pinning}, if the pinning condition \eqref{eqn:pinning-condition}, the KKT condition \eqref{eqn:KKT} and the second-order sufficient condition \eqref{eqn:2nd-test-pinning} hold at $p^*$, the corresponding framework is prestress stable.
	\end{proposition}

	\begin{proof}
		The proof comes directly by applying Lemmas \ref{lemma:pinning-kkt}, \ref{lemma:2nd-sufficient-pin} and Theorem \ref{thm:prestress-stability}.
	\end{proof}



	\begin{proof}[Proof of Lemma \ref{lemma:pinning-kkt}]
		Due to the complementary pinning condition \eqref{eqn:pinning-condition}, we can write $\nullspace \: R(p)$ as
		\begin{equation}\label{eqn:pinning-inf-flex}
			\nullspace \: R(p) = \nullspace \: R_g(p) + \mathcal{T}(p) = \nullspace \: R_t(p) \oplus \mathcal{T}(p),
		\end{equation}
		where the two matrices $R_g(p)$ and $R_t(p)$ are
		\begin{equation}\label{eqn:pinning-matrics}
			R_g(p)=\begin{pmatrix}
				R(p)\\
				G(p)
			\end{pmatrix}, \qquad R_t(p) = \begin{pmatrix}
				R(p)\\
				T(p)
			\end{pmatrix}.
		\end{equation}

		We observe that any left null-vectors of $R_t(p)$ must be the zero extension of the self-stresses of the rigidity matrix $R(p)$, i.e. for any $w \in \mathbb{R}^m, s \in \mathbb{R}^{D}$ with $w^T R(p) + s^T T(p) = 0$, we have $s=0$ and $w$ being a self-stress of $R(p)$. To see why, we observe that $s^T T(p) \in \mathcal{T}(p) \subset \nullspace \: R(p)$, while $s^T T(p) = -w^T R(p) \in \text{Range} (R(p)^T)$. Therefore, using Fredholm alternative, we have $s^T T(p) = w^T R(p) = 0$, meaning that $w$ is a self-stress.

		We now prove Lemma \ref{lemma:pinning-kkt} by a counting argument. Using \eqref{eqn:pinning-inf-flex}, we observe that $\dim \nullspace \: R_g(p)= \dim \nullspace \: R_t(p)$. Consequently, we also have $\dim \nullspace \: R_g^T(p) = \dim \nullspace \: R_t^T(p)$. Since any null vector of $R_t(p)$ is a zero extension of a self-stress and any zero extension of a self-stress is a left null vector of $ \nullspace \: R_g^T(p)$, the counting indicates that there are no additional left null vectors of $\nullspace \: R_g^T(p)$.
	\end{proof}

	\begin{proof}[Proof of Lemma \ref{lemma:2nd-sufficient-pin}]
		We observe that since the KKT condition \eqref{eqn:KKT} holds, the critical cones in the unpinned and pinned cases are the null spaces of $R_g(p^*)$ and $R_t(p^*)$ defined in \eqref{eqn:pinning-matrics}, i.e. $C_k(p^*) = \nullspace \: R_t(p^*)$ and $C_k^\text{pin}(p^*) = \nullspace \: R_g(p^*)$, where $C_k^\text{pin}(p^*)$ is defined in \eqref{eqn:critical-cone-pinning}. Therefore, using \eqref{eqn:pinning-inf-flex}, for any $v \in C_k^\text{pin}(p^*)$, we have $v = u + r$ with $u \in C_k(p^*)$ and $r \in \mathcal{T}(p^*)$. Since adding trivial flex does not change the value of the second-order test $v^T \nabla^2 L(p^*,\lambda^*) v$ (see \eqref{eqn:lagrangian-rigid-body}), we obtain the equivalence of the second-order sufficient conditions in the pinned and unpinned cases.

	\end{proof}

	\section{The numerical algorithm for the constrained optimization problem} \label{app:numerical-scheme}
	We provide our numerical algorithm for the constrained optimization problem \eqref{eqn:opt} and \eqref{eqn:opt-pinning}. Our algorithm uses a projected gradient descent method, which is the same one described in \cite{zappa2018monte}. To present our algorithm in a general form, we write our constrained optimization problem as
	\begin{align}
		\min_x & \qquad f(x), & & x \in \mathbb{R}^n,\\
		\text{s.t.} & \qquad c(x) = 0, & & c(x):\mathbb{R}^n \rightarrow \mathbb{R}^m,
	\end{align}
	where $f(x):\mathbb{R}^n \rightarrow \mathbb{R}$ is the objective function and $c(x)=0$ is the constrained set. Our constrained optimization algorithm can be described in \cref{alg:constrained-saddle-search}.
	\begin{algorithm}[H]
		\caption{Search for local minima of $f(x)$ with constraints $c(x)=0$}
		\label{alg:constrained-saddle-search}
		\begin{algorithmic}[1]
			\State \textbf{Input:} Initial value $x_0$, step size \( \eta \), and tolerance tol.
			\While{ $|x_{n} - x_{n-1}|>$ tol }
			\State Set $x_{n+1} = x_n - \eta \nabla f(x_n)$

			\State Update $x_{n+1}$ by projecting it to $c(x) = 0$ via Newton's method: solve for $\alpha_n \in \mathbb{R}^m$ such that $c(x_{n+1} + (\nabla c(x_n) )^T\alpha_n) = 0$ and then update $x_{n+1} = x_{n+1} + \nabla c(x_n) )^T\alpha_n$.

			\State Set $n = n+1$

			\EndWhile
		\end{algorithmic}
	\end{algorithm}

\section*{Acknowledgments}
The authors gratefully acknowledge valuable discussions and comments during the ICERM workshop Geometry of Materials (April 7-11, 2025), which have contributed to improving this work.

M.H.C. acknowledges support from the Natural Sciences and Engineering Research Council of Canada (NSERC), RGPIN-2023-04449 / Cette recherche a \'{e}t\'{e} financ\'{e}e par le Conseil de recherches en sciences naturelles et en g\'{e}nie du Canada (CRSNG). X.L. and C.D.S. receive no external funding for this project.

\section*{Declarations}
The original manuscript was completed by the end of 2025. In preparing this submission, the authors used OpenAI Codex to refine and document the MATLAB code, improve the presentation of figures, prepare the accompanying animations, and correct grammar. The authors assume responsibility for all content.

\bibliographystyle{plainurl}
\bibliography{ref}
\end{document}